%% file: Conference.tex
\documentclass{ifacconf}
\usepackage[dvipsnames]{xcolor}
\usepackage{natbib}
\usepackage{amsmath,amsfonts}
\usepackage{mathrsfs}
\usepackage{array}
\usepackage[font=scriptsize,labelfont=sf,textfont=sf]{subfig}
\usepackage{textcomp}
\usepackage{stfloats}
\usepackage{verbatim}
\usepackage{graphicx}
\usepackage{tikz}
\usepackage{balance}
\usepackage{bm}
    \usetikzlibrary{arrows}
    \usetikzlibrary{positioning}
\newcommand{\pfrac}[2]{\partial_{#2}{#1}} 
\newcommand{\diag}[1]{\text{diag}\{#1\}} 
\usepackage{caption,booktabs}

\newtheorem{thme}{Theorem}
\newtheorem{dfn}{Definition}
\newtheorem{prop}{Proposition}

\newtheorem{remm}{Remark}       
        \allowdisplaybreaks
\begin{document}
\begin{frontmatter}
\title{Output Feedback Control of Suspended Sediment Load  Entrainment in Water Canals and Reservoirs}
\thanks[footnoteinfo]{This work was funded by the NSF CAREER Award CMMI-2302030.}

\author[es]{Eranda Somathilake}
\author[md]{Mamadou Diagne}

\address[es]{University of California San Diego, CA 92093 USA (e-mail: esomathilake@ucsd.edu)}
\address[md]{University of California San Diego, CA 92093 USA (e-mail: mdiagne@ucsd.edu)}

\maketitle

\begin{abstract}
This paper addresses the management of water flow in a rectangular open channel, considering  the dynamic nature of both the channel's bathymetry and the suspended  sediment particles caused by  entrainment and deposition effects. The control-oriented model under study is a  set of   coupled nonlinear partial differential equations (PDEs) describing conservation of mass and  momentum while accounting for constitutive relations that govern sediment erosion and deposition phenomena.   The proposed boundary control problem presents a fresh perspective in water canal management and expands   Saint-Venant Exner (SVE) control frameworks by integrating dynamics related to the transport of fine particles. After linearization, PDE backstepping design is employed to stabilize both the bathymetry, the water dynamics together with the concentration of suspended sediment particles. Two underflow sluice gates are used for flow control at the upstream and downstream boundaries with only the downstream component being actuated. An observer-based backstepping control design is carried out for the downstream gate using state measurement at the upstream gate to globally exponentially stabilize the linearized system to a desired equilibrium point in $\mathscr{L}^2$ sense. The stability analysis is performed on the  linearized model which is a system of four coupled PDEs, three of which are rightward convecting and one leftward.  The proposed control design has the potential to facilitate efficient reservoir flushing operations. Consistent simulation results are presented to illustrate the feasibility of the designed control law.

\end{abstract}

\end{frontmatter}

Water infrastructures that obstruct natural sediment flow and disrupt natural watercourses are widely recognized as major contributors to the reduction of water reservoirs' storage capacity, primarily due to the accumulation of sediment particles. This poses a challenge in meeting the increasing demand for water  supply, especially in the face of the rapid and continuous growth of the world's population under the stress of climate change. Smart operation planning strategies of engineered rivers are needed to prevent fundamental transgenerational challenges such as full or partial dam decommissioning or  costly restoration of   reservoirs' storage capacity.
\section{Introduction}

Water, much like the regulated pulse of blood in arteries, navigates through the veins of a city via essential canals, necessitating a carefully designed flow control system. Civil engineering  infrastructures such as canals, are essential to the effective  distribution of water resources for industrial and domestic purposes, agriculture and energy generation. Sustainable management of water distribution systems ensures a consistent and  reliable water supply while prioritizing ecological preservation and low maintenance costs. Previous studies (\cite{pierre1998}) have considered the use of gate openings to regulate water flow in canals. The control design can be carried out such that the required discharge is achieved considering the dynamics of the flow in the canal. The dynamics of open-channel hydraulic systems obey the conservation of mass and momentum and consists of two first-order coupled partial differential equations (PDEs) (\cite{saint1871})  known as Saint-Venant equations. Numerous studies, such as \cite{balogun1988,coron1999,litrico2009,dos2012multi,nasir2021} (and the references therein) have extensively explored various control strategies for this PDE system, including model predictive control, internal model boundary control (IMBC), linear–quadratic regulator (LQR) and Lyapunov design. The initial result in \cite{coron1999} designed a controller for a system of coupled $2\times2$ hyperbolic PDEs, introducing an entropy-based Lyapunov function for stabilizing the Saint-Venant equations with two regulator gates. Subsequently, in \cite{bastin2010}, a controller design using a quadratic Lyapunov function under specific constraints was presented for a $2\times 2$ hyperbolic system. Both approaches required actuation at two boundaries for stabilization. In contrast, \cite{vazquez2011} introduced a boundary-stabilizing controller for a $2\times 2$ coupled linear hyperbolic PDE system in the canonical form, employing a collocated observer-based approach. However, these approaches, while successful in achieving stabilization, did not account for the dynamics of  sediment, both in \emph{suspended particle} form as well as \emph{bed load transport.} 

Sediment transportation primarily as bed load transport can be incorporated into the model by coupling the Saint-Venant equations with the Exner equation (\cite{Hudson2003}). Boundary control design of the Saint-Venant Exner (SVE) model has been explored in subsequent works. In \cite{ababacar2012} dissipative boundary conditions were assumed to guarantee $\mathscr{L}^2$-exponential stability of the closed-loop system whereas  \cite{tang2014} employ singular perturbation to conceive a  boundary controller  for hyperbolic systems with application to SVE. As the only result that enables stabilization under a supercritical flow regime with a single boundary actuation  and anti-collocated boundary sensing, \cite{ababacar2017} achieve observer-based exponential stability of an equilibrium point via the PDE backstepping approach. The result in \cite{ababacar2017} was followed by the stabilization of the bilayer model involving two separable fluids with different densities in \cite{mamadou2017}. The results in \cite{ababacar2017} and \cite{mamadou2017} are specific cases  of stabilizing systems with coupled  $n+1$ and $n+m$ heterodirectional linear hyperbolic PDEs (\cite{dimeglio2013} and \cite{hu2015control}).

The consideration of  sediment transport encompassing suspended sediment and bed load advection is crucial as discussed in \cite{owens2005} because sedimentation impacts the channel morphology, behaviour and navigation leading to heightened operational costs, diminished efficiency, and adverse effects on living species  habitats (\cite{ahn2008hydrology}). Recent discoveries revealed surface coarsening in gravel-bedded rivers when upstream bed load supply falls short of flow capacity  to transport the load (\cite{dietrich1989sediment}). In shaping river channel geometry, the relevancy of sediment motion thresholds is emphasized by \cite{phillips2022threshold}, while \cite{lamb2020mud} studied  the correlation between mud flocculation and suspended sediment particles. Since primarily, reservoir flushing occurs as suspended sediment (\cite{jihn1996}), a clear understanding of the behavior of sediment particles entrainment in flowing water is key to efficient flushing operations.  Controlled flushing operation   improves trap efficiency
in reservoirs by storing the clear water during low sediment concentration periods and releasing the turbid water
during the flood season  (\cite{idrees2019estimating}). The operational efficiency of an  automated sediment bypass tunnel, redirecting the inflow of fine sediment particles to Lake Miwa's 69 $m$ high gravity concrete dam and its 29.95 million cubic meters  reservoir capacity, is detailed in \cite{kantoush2011evaluation}. The model-free experimental results in \cite{kantoush2011evaluation} demonstrates the impact of sediment management automation in water infrastructures, allowing  for the bypassing of  399,000 $m^3$ of sediment annually, which thereby, extends  the life cycle of the subjacent multipurpose reservoir. Furthermore, a model based approach to reservoir flushing is proposed in \cite{jihn1996,tarekegn2014}, which eliminates the need for a complex monitoring system as a trade off for accuracy. Experimental validation of the model proposed by \cite{jihn1996} is carried out for a fixed gate opening and in \cite{tarekegn2014}, simulation results were shown to study the behavior of sediment deposition during reservoir flushing. The experimental and simulation studies described earlier are executed  in an open-loop system, which leaves the  feedback control design for suspended sediment flushing unaddressed.


More advanced representations of sediment transport account for  bed load transport, namely, the Exner equation, coupled with a $1D$ non-equilibrium \emph{suspended sediment transport equation}  to describe the sediment bed deformation dynamics (\cite{parker1986,garcia2008}). The  result is practically useful for conducting flushing operations of reservoirs.

We consider the coupled hyperbolic PDEs in \cite{parker1986,garcia2008} equipped with a model-based boundary controller that enables one to achieve a real-time regulation of  suspended sediment concentration and bed profile in a water canal. 
The equations governing the sediment dynamics illustrate how the motion of suspended sediment influences bed load deformation through entrainment and deposition phenomena. Mathematically, the model is governed by a system of $4\times4$ hyperbolic PDEs. The non-equilibrium suspended sediment transportation is governed by a $2\times2$ hyperbolic system of  water motion under a constant bottom slope combined with a $1D$  hyperbolic PDE reflecting the evolution of the suspended sediment concentration and another PDE for channel bed evolution. It is worth emphasizing that an experimental validation for the $3D$ version of \cite{parker1986,garcia2008} can be found  in \cite{zeng2005}. The nonlinear model in \cite{garcia2008} is linearized about a desired equilibrium point.  This point is determined by the required water demand, the optimal sediment concentration and a uniform channel bed with a constant slope.   The control design relies on an upstream and a downstream gate where the upstream gate is opened to supply the required discharge. The backstepping approach is carried out to determine the appropriate discharge on the downstream gate such that the system reaches the equilibrium point exponentially in $\mathscr{L}^2$ sense, which requires the design of an observer that exponentially converges in $\mathscr{L}^2$ sense as well.

\begin{figure}
\centering
\includegraphics[width=\columnwidth]{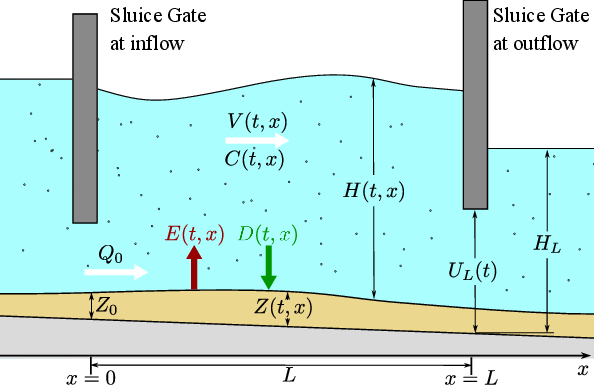}
\caption{Illustration of the system}\label{fig:sys}
\end{figure}

This article is organized as follows. The governing equations are introduced and linearized about the desired equilibrium point, then a coordinate transform is considered to facilitate the control design in section~\ref{sec:mathematical model}. The full state feedback control design is presented in section~\ref{sec:state feedback controller}. Next the observer design and the output feedback controller are given in sections~\ref{sec:backstepping obsever design} and \ref{sec:output feedback control}, respectively. The simulation results for the system with the output feedback controller are in section~\ref{sec:simulation results}. Finally, the concluding remarks are given in section~\ref{sec:conclusion}.

\begin{dfn}
\rm The $\mathscr{L}^2$-norm of a function $f:[0,L]\to\mathbb{R}$  is defined as
\begin{align*}
    \|f(x)\|_{\mathscr{L}^2}=\left(\int_0^L|f(x)|^2dx\right)^\frac{1}{2}.
\end{align*}
\end{dfn}

\section{Mathematical model}\label{sec:mathematical model}
We examine a morphodynamics model consisting of one-dimensional nonlinear hyperbolic PDEs describing unsteady flow, which considers suspended sediment of uniform grain size on a moving bed within a straight, prismatic open channel. This channel features a rectangular cross-section with a unit width and a length represented by $L$. The flow depth $H(t,x)$, the depth-averaged flow velocity $V(t,x)$, the depth averaged volumetric suspended sediment concentration $C(t,x)$, and the bathymetry $Z(t,x)$ which is the depth of the sediment layer above the channel bottom are the state variables of the system. Here $t$ is the time and $x$ is the longitudinal coordinate. The dynamics of the system is described by the coupling of the Saint-Venant equations, the suspended sediment transport equation, and the Exner equation for the conservation of the bed sediment. The bottom bed is assumed to undergo gradual variations. The equations are defined in \cite{garcia2008} and are characterized as
\begin{subequations}\label{PDE}
\begin{align}
\label{PDE a}
    \pfrac{H}{t}+V\pfrac{H}{x}+H\pfrac{V}{x}=&0,
    \\
\label{PDE b}
    \pfrac{V}{t}+g\pfrac{H}{x}+V\pfrac{V}{x}+g\pfrac{Z}{x}=&gS_b-C_f\frac{V^2}{H},
    \\
\label{PDE c}
    \pfrac{Z}{t}+ aV^2\pfrac{V}{x}=&\frac{\nu_s}{1-p^\prime}(D-E),
    \\
\label{PDE d}
    \pfrac{C}{t}+V\pfrac{C}{x}=&\frac{\nu_s}{H}(E-D),
\end{align}
\end{subequations}
where $g$ is the gravitational acceleration coefficient; $S_b$ is the constant slope of the channel bottom; $a$ is a constant given in Appendix \ref{sec:app_1}; $p^{\prime}$ is the sediment porosity ; $C_f$ is the friction coefficient assumed to be constant along the channel; $v_s$ is the sediment settlement velocity; $E(t,x)=E(V)$ and $D(t,x)=D(V,C)$ are the dimensionless sediment entrainment and deposition rates. Empirical formulas of $E$ and $D$ used in this study are (\cite{garcia2008})
\begin{align}\label{erosion}
    E=&\frac{A_1V^5}{1+A_2V^5},\\
    D=&\left(r_1+r_2\left(\frac{\sqrt{C_f}V}{v_s}\right)^{r_3}\right)C, \label{deposition}
\end{align}
here $A_1,A_2,r_1,r_2$, and $r_3$ are known constants that depend on the sediment and the water channel (see. Appendix \ref{sec:app_1}).

The upstream flow rate is assumed to have a constant value of $Q_0$ which can be maintained using the upstream gate. Furthermore, the channel bed elevation at the upstream gate can be assumed to be a constant value $Z_0$ which can be maintained by having an unerodable bed or the dredging operation (\cite{Cunge1980}). The suspended sediment concentration at the inlet is known which is modeled using sediment rating curves (\cite{Asselman2000}). Hence, the boundary conditions at $x=0$ are
\begin{subequations}\label{bc0}
\begin{align}
    \label{bc1}
    V(t,0)H(t,0)=&Q_0,\\
    \label{bc2}
    Z(t,0)=&Z_{0},\\
    \label{bc3}
    C(0,t)=&c_1+c_2(Q_0)^{c_3},
\end{align}
\end{subequations}
where $c_1,c_2$, and $c_3$ are coefficients of the suspended sediment rating curve \cite{Asselman2000}. The aperture of the underflow gate at the downstream, $U_L(t),$  is the control action, hence the boundary condition at $x=L$ is described by the sluice gate discharge equation (\cite{prabhata1992}) as
\begin{align}
    \nonumber
    &H(t,L)V(t,L)=\\
    \label{bc4}
    &(k_G\sqrt{2 g})(U_L(t)-Z(t,L))\sqrt{H(t,L)+Z(t,L)-H_L},
\end{align}    
where $k_G$ is a constant discharge coefficient and $H_L$ is the constant flow elevation beyond the downstream gate. This system is illustrated in Fig.~\ref{fig:sys}. The values of the parameters used in this study and their  relations are given in Appendix~\ref{sec:app_1}.


\subsection{Linearization and characteristic coordinates}
Note that the sediment entrainment and deposition rates will be equal at equilibrium (\cite{garcia2008}). We obtain constant values for equilibrium  states of the system defined as, $\begin{pmatrix}H_{eq}, V_{eq}, Z_{eq}, C_{eq}\end{pmatrix}^T$ which satisfy the relationship $gH_{eq}S_b=C_fV_{eq}^2$. Let the control input at the equilibrium be $U_{eq}$. Let the deviation of the states be $\mathcal{Y}=\begin{pmatrix}
        h,&v,&z,&c
    \end{pmatrix}^T=
    \begin{pmatrix}
        H-H_{eq},&V-V_{eq},&Z-Z_{eq},&C-C_{eq}
    \end{pmatrix}^T,$
and define the deviation of the control input from the equilibrium point as $u_L=U_L-U_{eq}$.
Then, the linearized system about the steady state is
\begin{align}
    \label{lin}
    \pfrac{\mathcal{Y}}{t}+\mathcal{A}\pfrac{\mathcal{Y}}{x}=\mathcal{B}\mathcal{Y},
\end{align}
where
\begin{align}
\nonumber
    \mathcal{A}=
    \begin{pmatrix}
        V_{eq}&H_{eq}&0&0\\
        g&V_{eq}&g&0\\
        0&aV_{eq}^2&0&0\\
        0&0&0&V_{eq}
    \end{pmatrix},\quad
    \mathcal{B}=
    \begin{pmatrix}
        0&0&0&0\\
        \phi_{vh}&\phi_{vv}&0&0\\
        0&\phi_{zv}&0&\phi_{zc}\\
        \phi_{ch}&\phi_{cu}&0&\phi_{cc}
    \end{pmatrix}.
\end{align}
The entries of the matrix   $\mathcal{B}$ are known constants that depend on the equilibrium (see Appendix~\ref{sec:app_2}). Similarly, the linearized boundary conditions are
\begin{subequations}\label{linbc}
\begin{align}
    \label{linbc1}
    \Pi\mathcal{Y}(t,0)=&\mathbf{0},\\
    \label{linbc2}
    \Pi_L\mathcal{Y}(t,L)+\pi_Lu_L(t)=&0,
\end{align}    
\end{subequations}
where
\begin{align}
    \nonumber
    \Pi=&\begin{pmatrix}
        V_{eq}&H_{eq}&0&0\\
        0&0&1&0\\
        -c_2c_3\frac{Q_0^{c_3}}{H_{eq}}&-c_2c_3\frac{Q_0^{c_3}}{V_{eq}}&0&1
    \end{pmatrix},\\
    \nonumber
    \Pi_L=&\begin{pmatrix}
        \pi_{Lh},&\pi_{Lv},&\pi_{Lz},&0
    \end{pmatrix}.
\end{align}
The entries in $\Pi_L$ are known constants that depend on the equilibrium  (given in Appendix~\ref{sec:app_2}).
Let the eigenvalues of $\mathcal{A}$ be $\lambda_k$ for $k=1,2,3,4$. For $k=1,2,3$, exact values can be obtained as shown in \cite{Hudson2003} and their approximate values for small $a$ is
\begin{align}
    \nonumber
    &\lambda_1 \approx V_{eq}-\sqrt{H_{eq}g},\quad\lambda_2 \approx \frac{agV_{eq}^3}{gH_{eq}-V_{eq}^2},\\
    &\lambda_3 \approx V_{eq}+\sqrt{H_{eq}g}.
\end{align}
The corresponding right eigenvectors, $R_k$ and approximate left eigenvectors, $L_k$ for $k=1,2,3$ can be calculated as
\begin{subequations}
\begin{align}
    R_k=&\begin{pmatrix}1,&\frac{\lambda_k-V_{eq}}{H_{eq}},&\frac{(V_{eq}-\lambda_k)^2-gH_{eq}}{gH_{eq}},&0\end{pmatrix}^T,\\
    \nonumber
    L_k=&\frac{1}{(\lambda_k-\lambda_i)(\lambda_k-\lambda_j)}\begin{pmatrix}(V_{eq}-\lambda_i)(V_{eq}-\lambda_j)+gH_{eq}\\H_{eq}\lambda_k\\gH_{eq}\\0\end{pmatrix}^T\\
    &\text{for }k\not=i\not=j,~i,j\in\{1,2,3\}.
\end{align}
\end{subequations}
The exact value for $\lambda_4$ and the corresponding right and left eigenvectors $R_4$ and $L_4$ respectively are
\begin{subequations}
\begin{align}
    &\lambda_4=V_{eq},\\
    &R_4=\begin{pmatrix}0,&0,&0,&1\end{pmatrix}^T,\quad L_4=\begin{pmatrix}0,&0,&0,&1\end{pmatrix}.
\end{align}    
\end{subequations}
Finally, defining 
\begin{align*}
\mathcal{R}=&\begin{pmatrix}R_1,&R_2,&R_3,&R_4\end{pmatrix},\\ \mathcal{L}=&\begin{pmatrix}L_1^T,&L_2^T,&L_3^T,&L_4^T\end{pmatrix}^T,
\end{align*} 
\eqref{lin} can be written in characteristic coordinates with $\Lambda=\diag{\lambda_1,\lambda_2,\lambda_3,\lambda_4}$ as
\begin{align}
    \label{characteristic}
    &\pfrac{\bm{\phi}}{t}+\Lambda\pfrac{\bm{\phi}}{x}=\mathcal{D}\bm{\phi},\\
\nonumber
    &\bm{\phi}=\mathcal{L}\mathcal{Y},\quad \mathcal{D}=\mathcal{L}\mathcal{B}\mathcal{R}, \quad \Lambda=\mathcal{L}\mathcal{A}\mathcal{R},
\end{align}
where $d_{ij}$ are the elements of row $i$ and column $j$ of $\mathcal{D}$.
\subsection{Control problem statement}
Using the system defined in \eqref{characteristic}, a control problem can be formulated which is similar to the problem solved in \cite{ababacar2017}. The behavior of the system at equilibrium (subcritical or supercritical flow) is characterized by the Froude number, $Fr$  defined as
\begin{align}
    \label{froude}
    Fr=\frac{V_{eq}}{\sqrt{gH_{eq}}}.
\end{align}
$\Lambda$ satisfy the inequality $\lambda_1<0<\lambda_2\ll\lambda_4<\lambda_3$ for subcritical flow($Fr<1$) and $\lambda_2<0<\lambda_1<\lambda_4<\lambda_3$ for supercritical flow($Fr>1$). 

\begin{remm}\rm Water flow within open channels exhibits two distinct behaviors: subcritical and supercritical flows. Subcritical flow is characterized by slower flow velocity and deeper water depth, while supercritical flow has a higher flow velocity and shallower depth. This study aims to propose control designs for both flow regimes. However, it is essential to note a fundamental characteristic of supercritical flow: disturbances occurring downstream do not propagate upstream, as emphasized in \cite{chanson2004} for a channel without a moving bathymetry. Although a moving bathymetry allows for control via a downstream gate, the practical feasibility of such control greatly depends on the properties of the moving bathymetry. Additionally, owing to the elevated flow velocity within the supercritical flow regime, erosion rates can be considerably high. Consequently, achieving stabilization of fluid flow within the supercritical regime poses substantial difficulty in most cases.
\end{remm}

Using the variable change $\bm{u_0}=\mathcal{T}_1\bm{\phi}$ for subcritical flow and $\bm{u_0}=\mathcal{T}_2\bm{\phi}$ for supercritical flow where $\mathcal{T}_1$ and $\mathcal{T}_2$ defined as
\begin{align}
    \mathcal{T}_1=\begin{pmatrix}
        0&1&0&0\\0&0&1&0\\0&0&0&1\\1&0&0&0\\
    \end{pmatrix},\quad
    \mathcal{T}_2=\begin{pmatrix}
        1&0&0&0\\0&0&1&0\\0&0&0&1\\0&1&0&0\\
    \end{pmatrix},
\end{align}
we arrive at the following system from \eqref{characteristic}
\begin{align}
\label{u sys}
    &\pfrac{\bm{u_0}}{t}+\bm{\gamma}\pfrac{\bm{u_0}}{x}=\bm{\sigma} \bm{u_0},\\
    &\bm{u_0}=\begin{pmatrix}u_1,u_2,u_3,u_4\end{pmatrix}^T,\quad\bm{\gamma}=\diag{\gamma_1,\gamma_2,\gamma_3,\gamma_4}.
\end{align}
  %
Define $\mu=-\gamma_4$ where $\mu>$ 0 and $\sigma_{ij}$ as the elements of row $i$ and column $j$ of $\bm{\sigma}$. Next, using the change of variable $w(t,x)=u_4(t,x)e^{\frac{\sigma_{44}}{\mu}x}$ and defining $\alpha_i(x)=\sigma_{i4}e^{-\frac{\sigma_{44}}{\mu}x}$ and $\theta_i(x)=\sigma_{4i}e^{\frac{\sigma_{44}}{\mu}x}$ for $i=1,2,3$, the system can be rewritten with the state variables $\begin{pmatrix}\bm{u},&w\end{pmatrix}^T=\begin{pmatrix}u_1,&u_2,&u_3,&w\end{pmatrix}^T$ for $i=1,2,3$ as
\begin{subequations}\label{sys}
\begin{align}
    \nonumber
    &\pfrac{u_i(t,x)}{t}+\gamma_i\pfrac{u_i(t,x)}{x}=\sum_{j=1}^3\sigma_{ij}u_j(t,x)\\
    \label{sys a}
    &\mkern200mu+\alpha_i(x)w(t,x),\\
    \label{sys d}
    &\pfrac{w(t,x)}{t}-\mu\pfrac{w(t,x)}{x}=\sum_{j=1}^3\theta_j(x)u_j(t,x),\\
    \label{bc 1}
    &u_i(t,0)=\delta_iw(t,0),\\
    \label{bc 2}
    &w(t,L)=\rho_1u_1(t,L)+\rho_2u_2(t,L)+U,
\end{align}
\end{subequations}
where the controller is designed considering $U=\pi_u u_L$ as the control input. The parameters $\delta_1$, $\delta_2$, $\delta_3$, $\rho_1$, $\rho_2$, and $\pi_u$ are constants to be calculated using the stated coordinate transform and boundary conditions \eqref{linbc}.
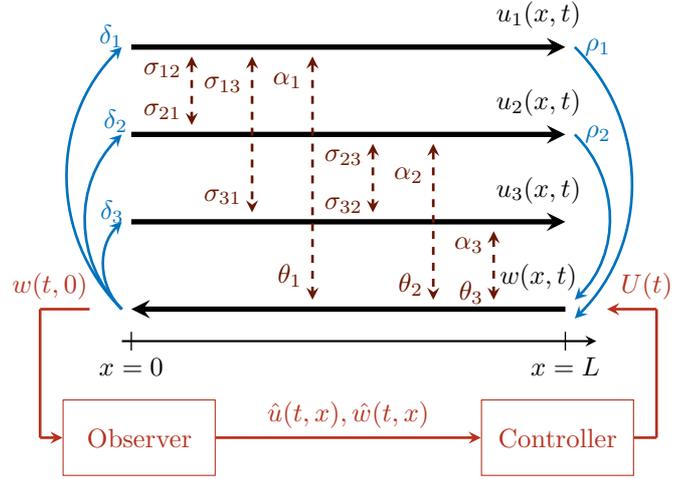
\begin{figure}[t]
    \centering
    \input{figures/system}
    \caption{Representation of the system to be controlled}
    \label{fig:system}
\end{figure}
\section{State feedback controller}\label{sec:state feedback controller}
Initially, we design a full-state feedback law which will then be implemented with the observed state variables.
\subsection{Backstepping transformation and the target system}
Consider the following backstepping transformation for $i=1,2,3$
\begin{subequations}\label{backstep}
\begin{align}
    \label{backstep1}
    \psi_i(t,x)=&u_i(t,x),\\
    \nonumber
    \chi(t,x)=&w(t,x)-\sum_{i=1}^3\int_0^xk_i(x,\xi)u_i(t,\xi)d\xi\\
    \label{backstep2}
    &-\int_0^xk_4(x,\xi)w(t,\xi)d\xi
\end{align}
\end{subequations}
to map the system \eqref{sys} to the following target system for $i=1,2,3$
\begin{subequations}\label{target}
\begin{align}
    \nonumber
    &\pfrac{\psi_i(t,x)}{t}+\gamma_1\pfrac{\psi_i(t,x)}{x}=\sum_{j=1}^3\sigma_{ij}\psi_j(t,x)\\
    \nonumber
    &\mkern20mu+\alpha_i(x)w(t,x)+\sum_{j=1}^3\int_0^xc_{ij}(x,\xi)\psi_j(t,\xi)d\xi\\
    \label{target a}
    &\mkern20mu+\int_0^x\kappa_i(x,\xi)\chi(t,\xi)d\xi,\\
    \label{target d}
    &\pfrac{\chi(t,x)}{t}-\mu\pfrac{\chi(t,x)}{x}= 0,\\
    \label{targetbc 1}
    &\psi_i(t,0)=\delta_i\chi(t,0),\\
    \label{targetbc 2}
    &\chi(t,L)=0,
\end{align}
\end{subequations}
where $c_{ij}(\cdot)$ and $\kappa_i(\cdot)$ are functions to be determined on the triangular domain $\mathbb{T}=\big\{(x,\xi)\in\mathbb{R}^2 | 0\leq\xi\leq x\leq L\big\}$. The sufficient condition for the backstepping transformation to hold can be found by substituting \eqref{backstep} to \eqref{target} and using \eqref{sys}, thereby it can be seen that the kernels $k_i(\cdot)$ should satisfy the following PDEs for $i=1,2,3$
\begin{subequations}\label{k}
\begin{align}
    \nonumber
    &\mu\pfrac{k_i(x,\xi)}{x}-\gamma_i\pfrac{k_i(x,\xi)}{\xi}=\sum_{j=1}^3\sigma_{ji}k_i(x,\xi)\\
    \label{k a}
    &\mkern220mu+\theta_i(\xi)k_4(x,\xi),\\
    \label{k d}
    &\mu\pfrac{k_4(x,\xi)}{x}+\mu\pfrac{k_4(x,\xi)}{\xi}=\sum_{i=1}^3\alpha_i(\xi)k_i(x,\xi),\\
    \label{k bc1}
    &k_i(x,x)=-\frac{\theta_i(x)}{\gamma_i+\mu},\\
    \label{k bc2}
    &\mu k_4(x,0)=\sum_{j=1}^3\gamma_j\delta_jk_j(x,0).
\end{align}
\end{subequations}
The existence, uniqueness, and continuity of the solutions for the kernel PDEs are shown in \cite{dimeglio2013}. The relations for $c_{ij}(\cdot)$ and $\kappa_i(\cdot)$ for all $i,j=1,2,3$ are defined as
\begin{align}
    \label{kappa}
    \kappa_i(x,\xi)=\alpha_i(x)k_4(x,\xi)+\int_\xi^x\kappa_i(x,s)k_4(s,\xi)ds,\\
    \label{cij}
    c_{ij}(x,\xi)=\alpha_i(x)k_j(x,\xi)+\int_\xi^x\kappa_i(x,s)k_j(s,\xi)ds.
\end{align}
The continuity of $k_i(\cdot)$ for $i=1,2,3,4$ implies the existence and continuity of the solution of the Volterra equations \eqref{kappa} as shown in \cite{linz1985}. Hence $c_{ij}(\cdot)$ and $\kappa_i(\cdot)$ for $i,j=1,2,3$ are continuous on $\mathbb{T}$ and thus bounded.

\subsection{Stability of the closed loop system}
To prove the stability of the closed loop system, we first show the exponential stability of the target system \eqref{target}. Consider the following Lyapunov function
\begin{align}
    \label{lyapunov}
    V_1(t)&=\int_0^L\left[p_1^ve^{-q_1^vx}\sum_{i=1}^3\frac{\psi_i(t,x)^2}{\gamma_i}+\frac{1+x}{\mu}\chi(t,x)^2\right]dx,
\end{align}
where $p_1^v$ and $q_1^v$ are positive parameters selected such that $\dot{V}_1(t)\leq 0$. Therefore, it can be concluded that the system \eqref{target} is exponentially stable.

\subsection{Inverse transformation and control law}
To ensure the equivalence of stability between the systems \eqref{sys} and \eqref{target}, the transformation \eqref{backstep} should be invertible. The invertibility of \eqref{backstep1} is trivial, using \eqref{backstep1}, transformation \eqref{backstep2} can be rewritten as
\begin{align}
\label{inv1}
    w(t,x)-\int_0^xk_4(x,\xi)w(t,\xi)d\xi=\Gamma(t,x),
\end{align}
with $\Gamma(t,x)=\chi(t,x)+\sum_{i=1}^3\int_0^xk_i(x,\xi)\psi_i(t,\xi)d\xi$. Now let
\begin{align}
\label{inv2}
    w(t,x)=\Gamma(t,x)+\int_0^xl_4(x,\xi)\Gamma(t,\xi)d\xi,
\end{align}
where $l_4(\cdot)$ can be determined as
\begin{align}
    l_4(x,\xi)=k_4(x,\xi)+\int_\xi^xl_4(s,\xi)k_4(x,s)ds.\label{l_4}
\end{align}
Next, from \eqref{inv2} substituting the value of $\Gamma(t,x)$, the inverse transformation can be determined as
\begin{align}
     w(t,x)=&\sum_{i=1}^3\int_0^xl_i(x,\xi)\psi_i(t,\xi)d\xi+\int_0^xl_4(x,\xi)\chi(t,\xi)d\xi,
 \end{align}
 where, for $i=1,2,3$
  \begin{align}
     l_i(x,\xi)=k_i(x,\xi)+\int_\xi^x l_4(x,s)k_i(s,\xi)ds.
 \end{align}
 The continuity of $k_i(\cdot)$ for $i=1,2,3,4$ implies the existence and continuity of the solution of the Volterra equation \eqref{l_4} as shown in \cite{linz1985}. Hence $l_i(\cdot)$ for $i=1,2,3,4$ are continuous on $\mathbb{T}$ and thus bounded.
Using \eqref{bc 2} and \eqref{targetbc 2} and the transformation \eqref{backstep}, we obtain the following control input 
\begin{align}
    \nonumber
    U(t)=&\sum_{i=1}^3\int_0^Lk_i(L,\xi)u_i(t,\xi)d\xi+\int_0^Lk_4(L,\xi)w(t,\xi)d\xi\\
    \label{U}
    &-\rho_1u_1(t,L)-\rho_2u_2(t,L).
 \end{align}
 that exponentially stabilizes the system \eqref{sys} in the $\mathscr{L}^2$ sense.
 Hence, the following theorem holds.
 \begin{thme}
 The system \eqref{sys} under the control law \eqref{U}, for any given initial condition ~~$\begin{pmatrix}\bm{u}(0,x),w(0,x)\end{pmatrix}~~\in(\mathscr{L}^2[0,L])^4$, is exponentially stable in the $\mathscr{L}^2$ sense.
 \end{thme}

\section{Backstepping observer design}\label{sec:backstepping obsever design}
The feedback controller \eqref{U} requires full state information along the spatial domain, but we consider that only the boundary value $y(t)=w(t,0)$ is available. Consider the following observer for $i=1,2,3$
\begin{subequations}\label{obs}
\begin{align}
    \nonumber
    &\pfrac{\hat{u}_i(t,x)}{t}+\gamma_i\pfrac{\hat{u}_i(t,x)}{x}=\sum_{j=1}\sigma_{ij}\hat{u}_j(t,x)+\alpha_i(x)\hat{w(t,x)}\\
    \label{obs a}
    &\mkern180mu-p_i(x)(y(t)-\hat{w}(t,0)),\\
    \nonumber
    &\pfrac{\hat{w}(t,x)}{t}-\mu\pfrac{\hat{w}(t,x) }{x}=\sum_{j=1}^3\theta_j(x)\hat{u}_j\\
    \label{obs d}
    &\mkern180mu-p_4(x)(y(t)-\hat{w}(t,0)),\\
    \label{obs bc1}
    &\hat{u}_i(t,0)=\delta_iy(t),\\
    \label{obs bc2}
    &\hat{w}(t,L)=\rho_1\hat{u}_1(t,L)+\rho_2\hat{u}_2(t,L)+U(t),
\end{align}
\end{subequations}
where $\begin{pmatrix}\bm{\hat{u}}&\hat{w}\end{pmatrix}^T=\begin{pmatrix}\hat{u}_1&\hat{u}_2&\hat{u}_3&\hat{w}\end{pmatrix}^T$ is the estimated state vector. The observer gains,  $p_1(x)$, $p_2(x)$, $p_3(x)$, and $p_4(x)$ are functions to be determined such that$\begin{pmatrix}\bm{\hat{u}}&\hat{w}\end{pmatrix}^T$ converges to $\begin{pmatrix}\bm{u}&w\end{pmatrix}^T$. Defining the error term as $\begin{pmatrix}\bm{\tilde{u}}&\tilde{w}\end{pmatrix}^T=\begin{pmatrix}\tilde{u}_1&\tilde{u}_2&\tilde{u}_3&\tilde{w}\end{pmatrix}^T$ as $\begin{pmatrix}\bm{\tilde{u}}&\tilde{w}\end{pmatrix}^T=\begin{pmatrix}\bm{\hat{u}}&\hat{w}\end{pmatrix}^T-\begin{pmatrix}\bm{u}&w\end{pmatrix}^T$, we obtain the following error dynamics for $i=1,2,3$
\begin{subequations}\label{obs err}
\begin{align}
    \nonumber
    &\pfrac{\tilde{u}_i(t,x)}{t}+\gamma_i\pfrac{\tilde{u}_i(t,x)}{x}=\sum_{j=1}^3\sigma_{ij}\tilde{u}_j(t,x)+\alpha_i(x)\tilde{w}(t,x)\\
    \label{obs err a}
    &\mkern200mu+p_i(x)\tilde{w}(t,0),\\
    \nonumber
    &\pfrac{\tilde{w}(t,x)}{t}-\mu\pfrac{\tilde{w}(t,x)}{x}=\sum_{j=1}^3\theta_j(x)\tilde{u}_j(t,x)\\
    \label{obs err d}
    &\mkern60mu+p_4(x)\tilde{w}(t,0),\\
    \label{obs err bc1}
    &\tilde{u}_i(t,0)=0,\\
    \label{obs err bc2}
    &\tilde{w}(t,L)=\rho_1\tilde{u}_1(t,L)+\rho_2\tilde{u}_2(t,L).
\end{align}
\end{subequations}
An invertible backstepping transformation is used to transform the system to an exponentially stable target system.The observer gains are determined such that the two systems map to each other.
\subsection{Backstepping transformation and the target system}
Consider the following backstepping transformation for $i=1,2,3$
\begin{subequations}\label{obs backstep}
\begin{align}
    \label{obs backstep1}
    \tilde{u}_i(t,x)=&\tilde{\psi}_i(t,x)+\int_0^xm_i(x,\xi)\tilde{\chi}(t,\xi)d\xi,\\
    \label{obs backstep2}
    \tilde{w}(t,x)=&\tilde{\chi}(t,x)+\int_0^xm_4(x,\xi)\tilde{\chi}(t,\xi)d\xi
\end{align}    
\end{subequations}
to map the system \eqref{obs err} to the following target system for $i=1,2,3$
\begin{subequations}\label{obs target}
\begin{align}
    \nonumber
    &\pfrac{\tilde{\psi}_i(t,x)}{t}+\gamma_i\pfrac{\tilde{\psi}_i(t,x)}{x}=\sum_{j=1}^3\sigma_{ij}\tilde{\psi}_j(t,x)\\
    \label{obs target a}
    &\mkern150mu+\sum_{j=1}^3\int_0^xg_{ij}(x,\xi)\tilde{\psi}_j(t,\xi)d\xi,\\
    \nonumber
    &\pfrac{\tilde{\chi}(t,x)}{t}-\mu\pfrac{\tilde{\chi}(t,x)}{x}=\sum_{j=1}^3\theta_j(x)\tilde{\psi}_j(t,x)\\
    \label{obs target d}
    &\mkern20mu\sum_{j=1}^3\theta_j(x)\tilde{\psi}_j+\int_0^xh_j(x,\xi)\tilde{\psi}_j(t,\xi)d\xi,\\
    \label{obs targetbc1}
    &\tilde{\psi}_i(t,0)=0,\\
    \label{obs targetbc2}
    &\tilde{\chi}(t,L)=\rho_1\tilde{\psi}_1(t,L)+\rho_2\tilde{\psi}_2(t,L),
\end{align}
\end{subequations}
where $g_{ij}(\cdot)$ and $h_{ij}(\cdot)$ are functions to be determined on $\mathbb{T}$. The sufficient conditions for the backstepping transformation to hold can be found by substituting \eqref{obs backstep} to the error system \eqref{obs err} and using the target system \eqref{obs target}, thereby it can be seen that the kernels $m_i(\cdot)$ should satisfy the following PDEs for $i=1,2,3$
\begin{subequations}\label{m}
\begin{align}
    \nonumber
    &\gamma_i\pfrac{m_i(x,\xi)}{x}-\mu\pfrac{m_i(x,\xi)}{\xi}=\sum_{j=1}^3\sigma_{ij}m_j(x,\xi)\\
    &\mkern230mu+\alpha_i(x)m_4(x,\xi),\\
    &\mu\pfrac{m_4(x,\xi)}{x}+\mu\pfrac{m_4(x,\xi)}{\xi}=-\sum_{j=1}^3\theta_j(x)m_j(x,\xi),\\
    &m_i(x,x)=\frac{\alpha_i(x)}{\gamma_i+\mu},\\
    &m_4(L,\xi)=\rho_1m_1(L,\xi)+\rho_2m_2(L,\xi).
\end{align}
\end{subequations}
    The existence, uniqueness, and continuity of the solutions for the kernel PDEs can be seen after the coordinate change $x=L-\bar{x},\xi=L-\bar{\xi}$, where we obtain a system similar to \eqref{k}. Also, the observer gains are
\begin{align}
    \label{obs gains}
    p_i(x)=-\mu m_i(x,0),~\text{for}~i=1,2,3,4.
\end{align}
The relationship for $g_{ij}(\cdot)$ and $h_i(\cdot)$ on $\mathbb{T}$ for $i,j=1,2,3$ is
\begin{subequations}
\begin{align}
    \label{h_i}
    h_i(x,\xi)=&-\theta_i(\xi)m_4(x,\xi)-\int_\xi^xm_4(x,s)h_i(s,\xi)ds,\\
    \label{g_ij}
    g_{ij}(x,\xi)=&-\theta_j(\xi)m_i(x,\xi)-\int_\xi^xm_i(x,s)h_j(s,\xi)ds.
\end{align}
\end{subequations}
The continuity of $m_i(\cdot)$ for $i=1,2,3,4$ implies the existence and continuity of the solution of the Volterra equations \eqref{h_i} as shown in \cite{linz1985}. Hence $g_{ij}(\cdot)$ and $h_i(\cdot)$ for $i,j=1,2,3$ are continuous on $\mathbb{T}$ and thus bounded.
\subsection{Exponential convergence of the observer}
Consider the following Lyapunov function
\begin{align}
    \label{obs lyapunov}
    V_2(t)&=\int_0^L\left[p_2^ve^{-q_2^vx}\sum_{i=1}^3\frac{\tilde{\psi}_i(t,x)^2}{\gamma_i}+\frac{e^{q_2^vx}}{\mu}\tilde{\chi}(t,x)^2\right]dx,
\end{align}
where $p_2^v,q_2^v>0$ are selected such that $\dot{V}_2\leq 0 $.Hence the target system \eqref{obs target} is exponentially stable.

\subsection{Inverse transformation}
The exponential stability of the origin of the error system \eqref{obs err} ensures that the states of the observer system \eqref{obs} converges to the system states exponentially. Since, the target system \eqref{obs target} is proven to be exponentially stable, the invertibility of the transformation \eqref{obs backstep} ensures that the observer states converge to the system states exponentially provided that the observer gains are defined as in \eqref{obs gains}. Let the inverse transformation for $i=1,2,3$ be
\begin{subequations}\label{obs inv}
\begin{align}
    \label{obs inv1}
    \tilde{\psi}_i(t,x)=&\tilde{u}_i(t,x)+\int_0^xn_i(x,\xi)\tilde{w}(t,\xi)d\xi,\\
    \label{obs inv2}
    \tilde{\chi}(t,x)=&\tilde{w}(t,x)+\int_0^xn_4(x,\xi)\tilde{w}(t,\xi)d\xi,
\end{align}
\end{subequations}
where the inverse kernels $n_i(\cdot)~\text{for}~i=1,2,3,4$ are defined as
\begin{subequations}\label{n}
\begin{align}
    n_4(x,\xi)=&m_4(x,\xi)+\int_\xi^xm_4(x,s)n_4(s,\xi)ds,\label{n_4}\\
    n_i(x,\xi)=&-m_i(x,\xi)-\int_\xi^xm_i(x,s)n_4(s,\xi)ds,~i\neq4.\label{n_i}
\end{align}
\end{subequations}
The continuity of $m_i(\cdot)$ for $i=1,2,3,4$ implies the existence and continuity of the solution of the Volterra equation \eqref{n_4} as shown in \cite{linz1985}. Hence $n_i(\cdot)$ for $i=1,2,3,4$ are continuous on $\mathbb{T}$ and thus bounded. 
Hence the following theorem holds.
\begin{thme}
    The observer system \eqref{obs}, with the coefficient functions $p_i(x)$ for $i=1,2,3,4$ as defined in \eqref{obs gains}, for any given initial condition $\begin{pmatrix}\bm{\hat{u}}(0,x),\hat{w}(0,x)\end{pmatrix}\in(\mathscr{L}^2[0,L])^4$, exponentially converges to the system \eqref{sys} in the $\mathscr{L}^2$ sense.
\end{thme}

\section{Output feedback control}\label{sec:output feedback control}
The controller, \eqref{U} which requires full state information is implemented with the observer \eqref{obs} which estimates the real system states of the system \eqref{sys} using the measured output $y(t)$. Therefore, we design an output feedback controller. Consider the combined system with the states $\begin{pmatrix}\bm{u}&w&\bm{\hat{u}}&\hat{w}\end{pmatrix}^T$, this is equivalent to the system with states $\begin{pmatrix}\bm{\hat{u}}&\hat{w}&\bm{\tilde{u}}&\tilde{w}\end{pmatrix}^T$. Consider the following backstepping transformation for the states $\begin{pmatrix}\bm{\hat{u}}&\hat{w}\end{pmatrix}^T$ for $i=1,2,3$
\begin{subequations}\label{backstephat}
\begin{align}
    \label{backstephat1}
    \hat{\psi}_i(t,x)=&\hat{u}_i(t,x),\\
    \nonumber
    \hat{\chi}(t,x)=&\hat{w}(t,x)-\sum_{i=1}^3\int_0^xk_i(x,\xi)\hat{u}_i(t,\xi)d\xi\\
    \label{backstephat2}
    &-\int_0^xk_4(x,\xi\hat{w}(t,\xi)d\xi
\end{align}
\end{subequations}
to map the system \eqref{obs} to the following target system for $i=1,2,3$
\begin{subequations}\label{oftarget}
\begin{align}
    \nonumber
    &\pfrac{\hat{\psi}_i(t,x)}{t}+\gamma_1\pfrac{\hat{\psi}_i(t,x)}{x}=\sum_{j=1}^3\sigma_{ij}\hat{\psi}_j(t,x)\\
    \nonumber
    &\mkern0mu+\alpha_i(x)\hat{w}(t,x)+\sum_{j=1}^3\int_0^xc_{ij}(x,\xi)\hat{\psi}_j(t,\xi)d\xi\\
    \label{oftarget a}
    &+\int_0^x\kappa_i(x,\xi)\hat{\chi}(t,\xi)d\xi-p_i(x)\tilde{w}(t,0),\\
    \nonumber
    &\pfrac{\hat{\chi}(t,x)}{t}-\mu\pfrac{\hat{\chi}(t,x)}{x}=\\
    \label{oftarget d}
    &\left(-p_4(x)+\int_0^x\sum_{i=1}^3k_i(x,\xi)p_i(\xi)d\xi\right)\tilde{w}(t,0),
\end{align}
\end{subequations}
where the kernels $k_i(\cdot)$ for $i=1,2,3,4$ are  defined in \eqref{backstep}. Therefore, the transformation is invertible and satisfies the relations in \eqref{k}. Next, using the backstepping transformation \eqref{obs backstep}, the system \eqref{obs err} can be mapped to the target system \eqref{obs target} as shown above. Now, for the system with the combined states $\hat{\psi}_i(t,x),\hat{\chi}(t,x),\tilde{\psi}_i(t,x),$ and $\tilde{\chi}(t,x)$ for $i=1,2,3$, consider the following Lyapunov function
\begin{align}
    \nonumber
    V(t)=&\int_0^L\left[\hat{p}e^{-qx}\sum_{i=1}^3\frac{\hat{\psi}_i(t,x)^2}{\gamma_i}dx+\frac{e^{qx}}{\mu}\hat{\chi}(t,x)^2\right]dx\\
    \label{lyapunov_main}
    &+\int_0^L\left[\tilde{p}e^{-qx}\sum_{i=1}^3\frac{\tilde{\psi}_i(t,x)^2}{\gamma_i}dx+\frac{e^{qx}}{\mu}\tilde{\chi}(t,x)^2\right]dx,
\end{align}
where $\hat{p},\tilde{p},q>0$ are constants that are selected such that $\dot{V}(t)\leq0$. Hence we arrive at the following theorem.
\begin{thme}
The combined observer-model system with states $\begin{pmatrix}\bm{u}&w&\bm{\hat{u}}&\hat{w}\end{pmatrix}^T$ and with initial conditions 
\begin{align*}\begin{pmatrix}\bm{u}(0,x)&w(0,x)&\bm{\hat{u}}(0,x)&\hat{w}(0,x)\end{pmatrix}^T\in(\mathscr{L}^2[0,L])^8,\end{align*}
with the control law
    \begin{align}
        \nonumber
        U(t)=&\sum_{i=1}^3\int_0^Lk_i(L,\xi)\hat{u}_i(t,\xi)d\xi+\int_0^Lk_4(L,\xi)\hat{w}(t,\xi)d\xi\\
        \label{U of}
        &-\rho_1\hat{u}_1(t,L)-\rho_2\hat{u}_2(t,L),
    \end{align}
is exponentially stable in the $\mathscr{L}^2$ sense.
\end{thme}
\section{Simulation results}\label{sec:simulation results}
The simulations are performed for the subcritical flow regime, subject to the  output-feedback control \eqref{U of}. Consider a channel of unit length ($L=1$). The equilibrium states considered are: $H_{eq}=2~m,V_{eq}=3~m/s,Z_{eq}=0.4~m,$ and $C_{eq}=0.005$. Using \eqref{Ueq}, and the equilibrium values of the state variables, we obtain the gate aperture at equilibrium as $U_{eq}=2.308~m$. The initial bottom topography of the channel is defined as
\begin{align}\label{ini1}
    Z(0,x)=Z_{eq}+0.2\left(e^{-\frac{(x-L/2)^2}{0.1}}-e^{\frac{L^2}{0.4}}\right)
\end{align}
and the initial water height, water velocity, and suspended sediment concentration variation are defined as
\begin{align}\label{ini2}
    H(0,x)=&0.1\sin{\left(\frac{\pi x}{L}\right)}+H_{eq}+Z_{eq}-Z(0,x),\\
    V(0,x)=&\frac{Q_0}{H(0,x)}\label{ini3},\\
    C(0,x)=&0.003\sin{\left(\frac{\pi x}{2L}\right)} + C_{eq},\label{ini4}
\end{align}
respectively. The initial states of the observer are taken assuming that the system is at equilibrium. The kernel PDEs \eqref{k} and \eqref{m} are solved with the method presented in \cite[Appendix]{anfinsen2019adaptive}.

\begin{figure}
\centering
\subfloat[Norm of the characteristics for the closed loop system \label{fig:norm_cl}]{
\includegraphics[width=0.45\columnwidth,clip]{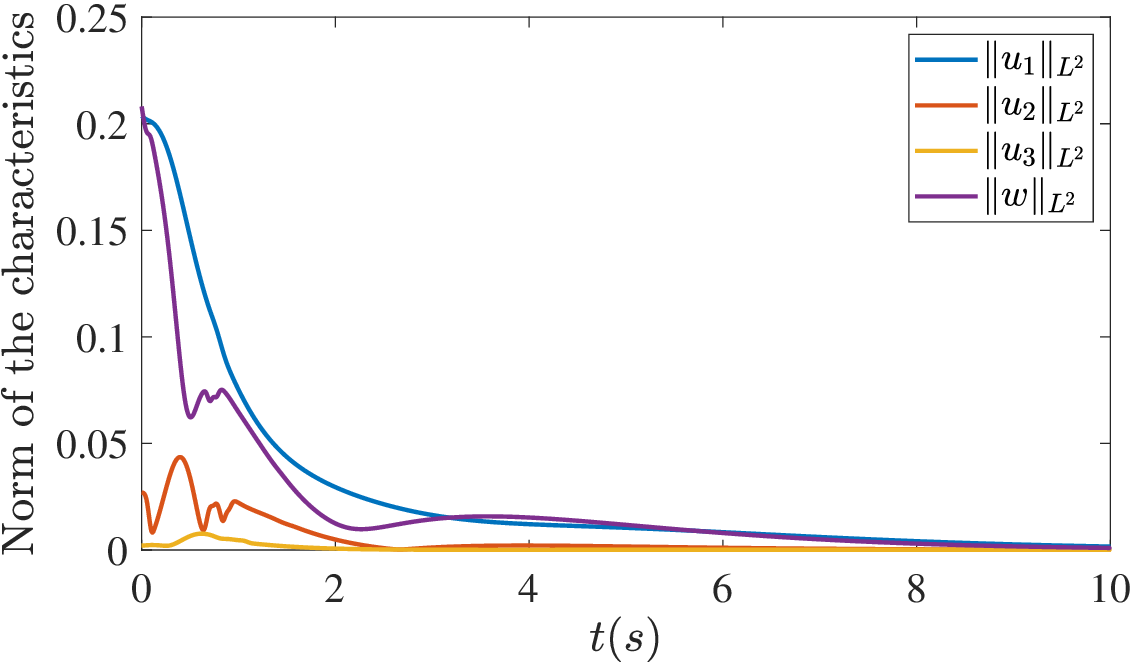}}
~
\subfloat[Norm of the characteristics for the open loop system \label{fig:norm_ol}]{
\includegraphics[width=0.45\columnwidth,clip]{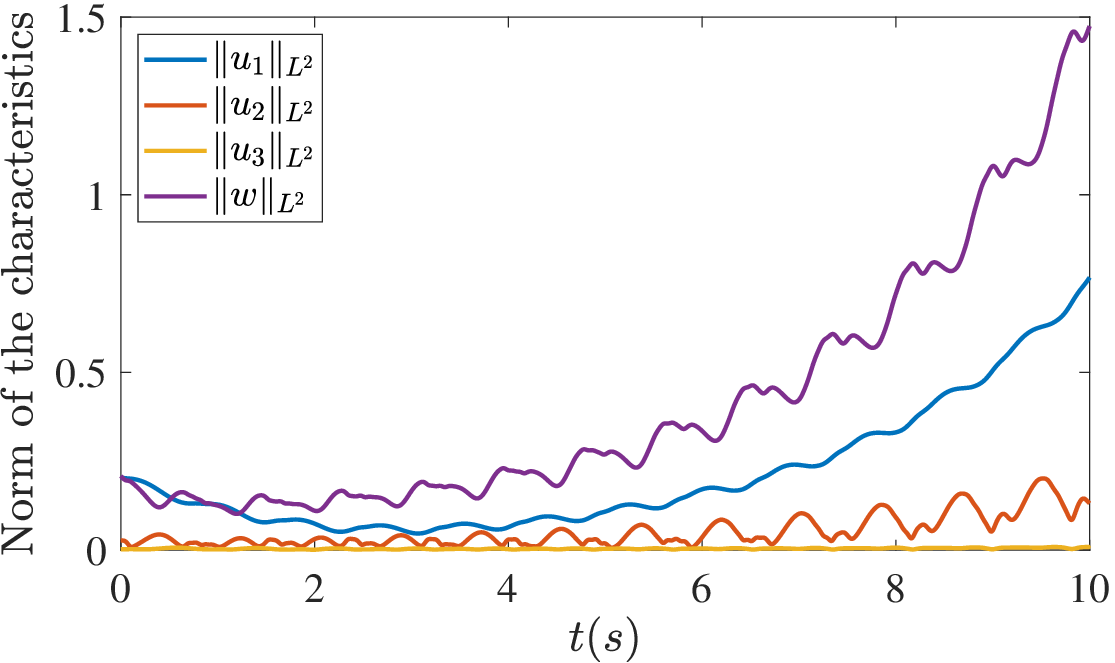}}
\caption{$\mathscr{L}^2$-Norm of characteristics for the open and closed loop systems}
\end{figure}

\begin{figure}
\centering
\subfloat[The control input $U_L(t)$ and the water level at the actuation boundary \label{fig:Control}]{
\includegraphics[width=0.45\columnwidth,clip]{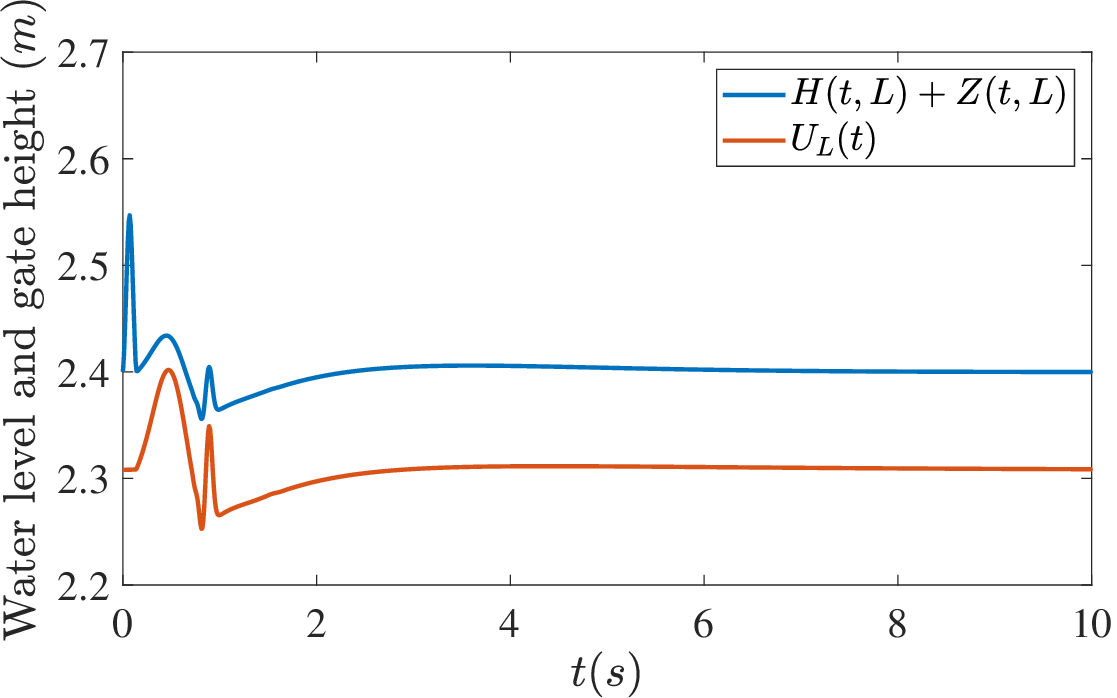}}
~
\subfloat[Water velocity at the actuation boundary \label{fig:v_L}]{
\includegraphics[width=0.45\columnwidth,clip]{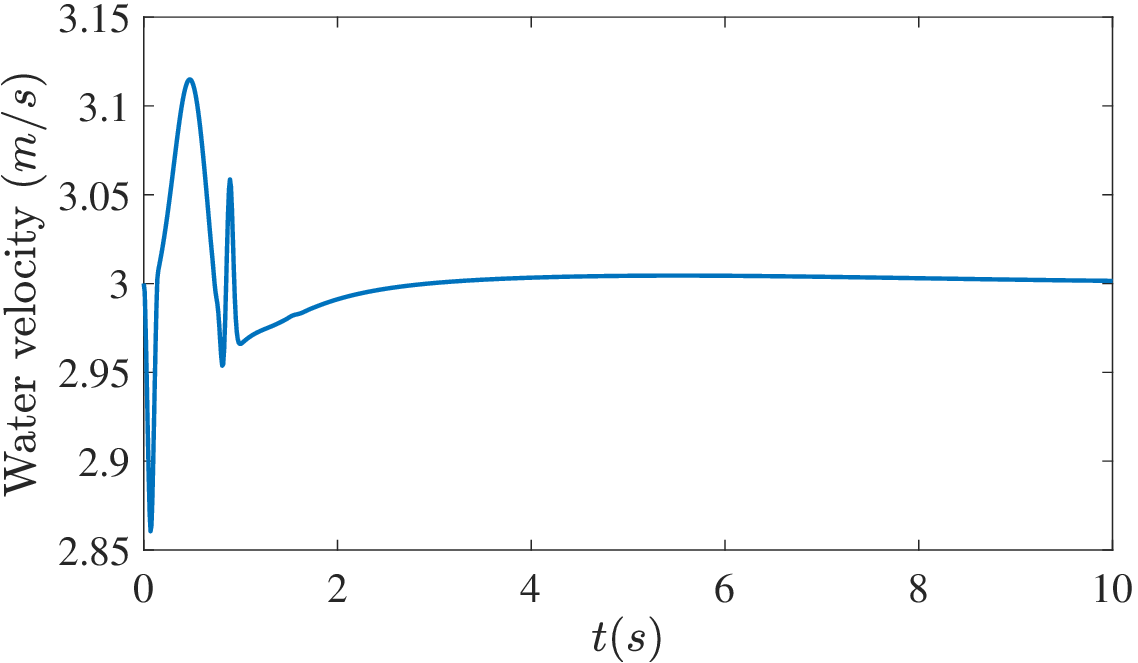}}
\caption{Water level, water velocity at the actuation boundary and the downstream gate height, $U_L(t)$}
\end{figure}

\begin{figure}
\centering
\subfloat[Dynamics of the water height \label{fig:h}]{
\includegraphics[width=0.45\columnwidth,trim={0 0 0 2cm},clip]{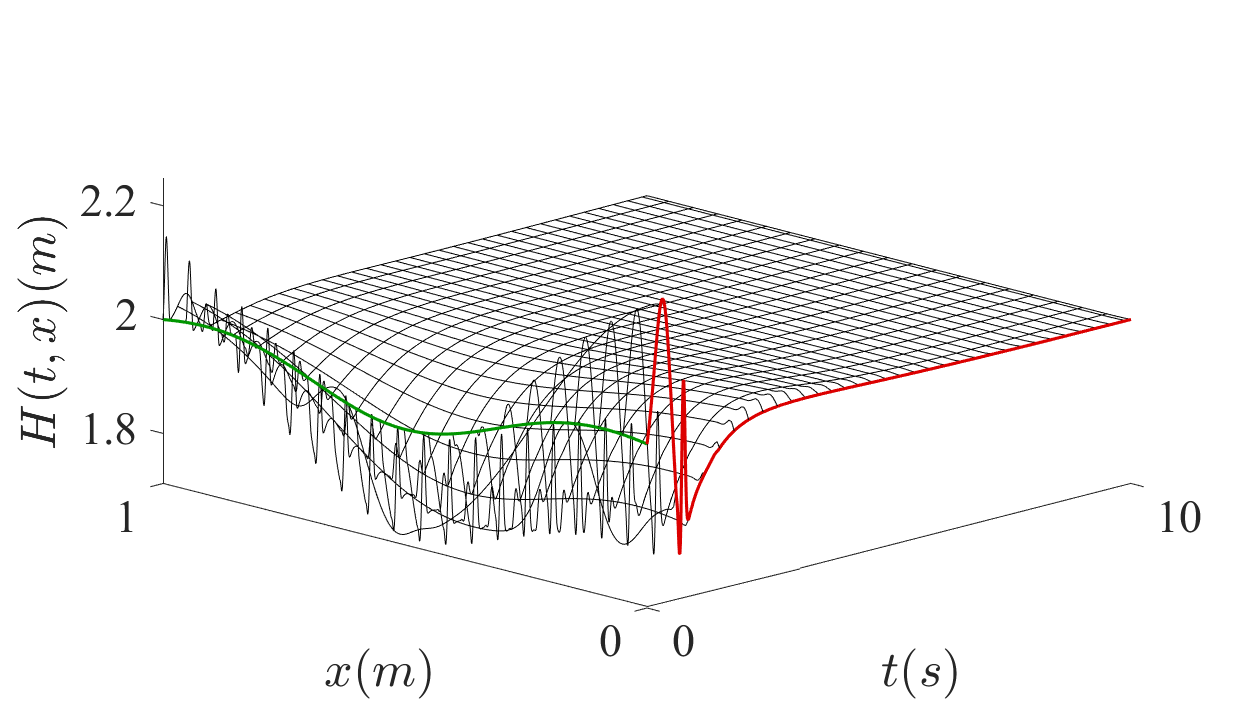}}
~
\subfloat[Dynamics of the water velocity \label{fig:v}]{
\includegraphics[width=0.45\columnwidth,trim={0 0 0 2cm},clip]{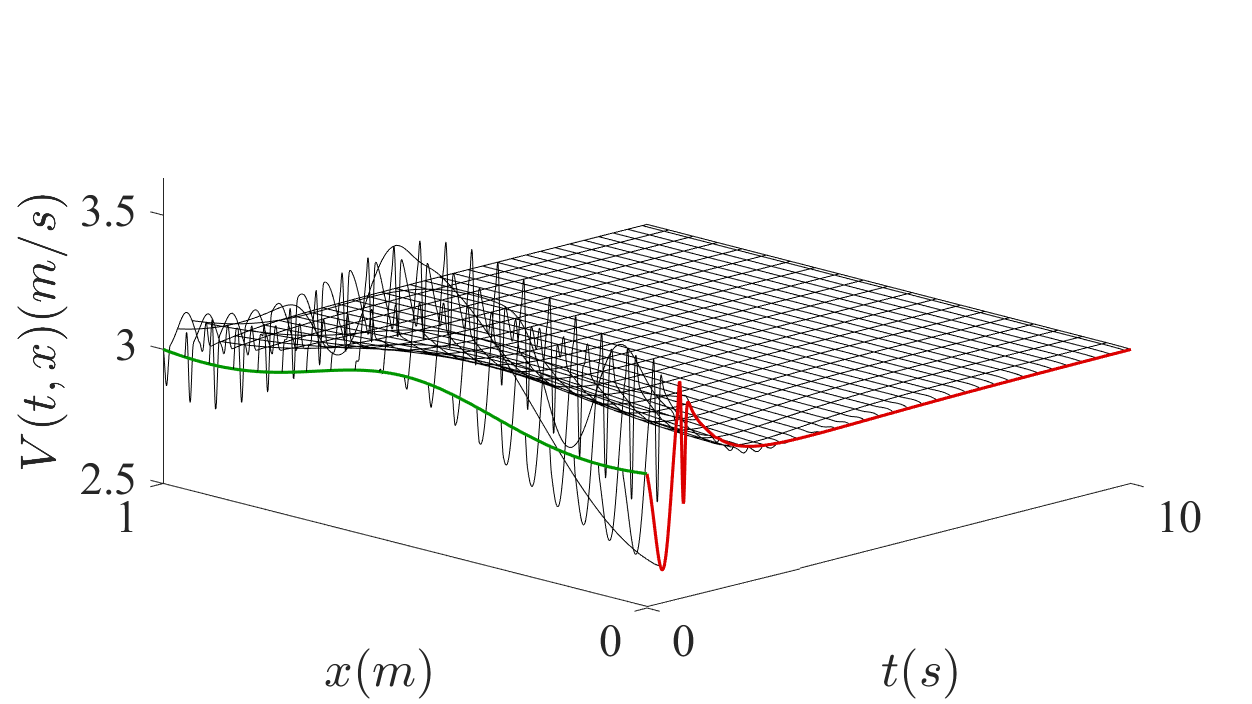}}
\\
\subfloat[Dynamics of the bathymetry   \label{fig:z}]{
\includegraphics[width=0.45\columnwidth,trim={0 0 0 2cm},clip]{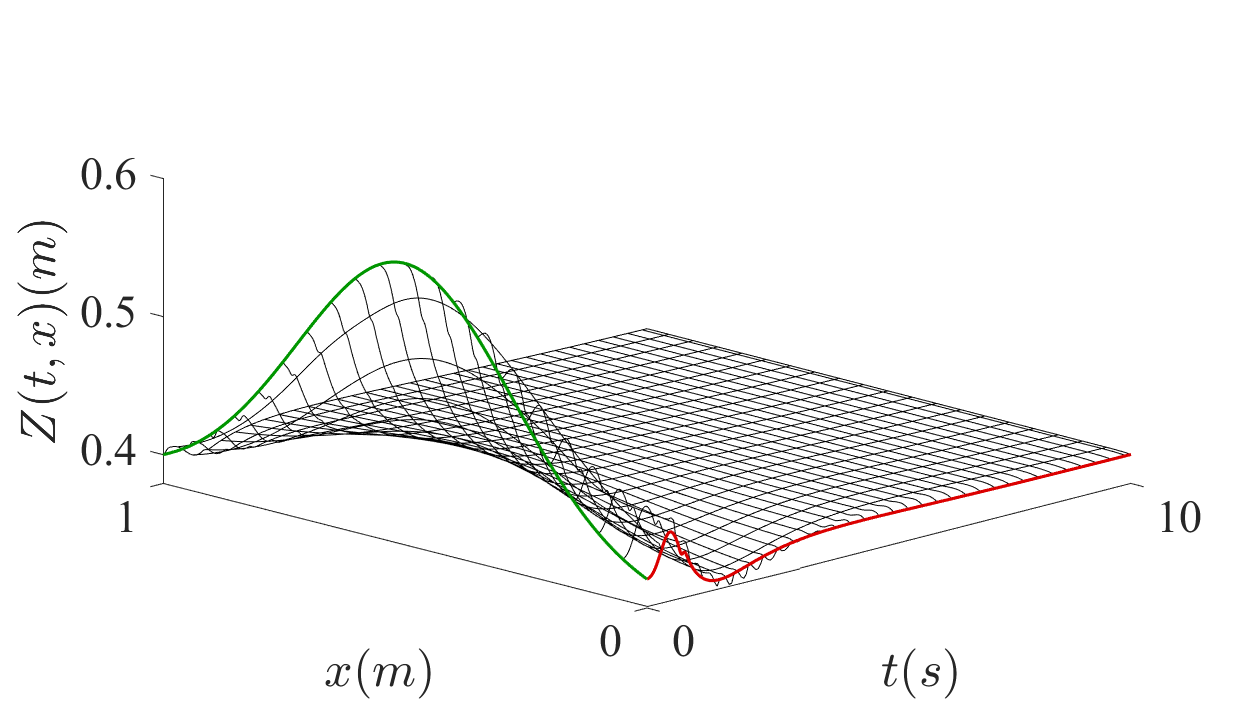}}
~
\subfloat[Dynamics of the concentration of suspended sediment particles \label{fig:c}]{
\includegraphics[width=0.45\columnwidth,trim={0 0 0 0.5cm},clip]{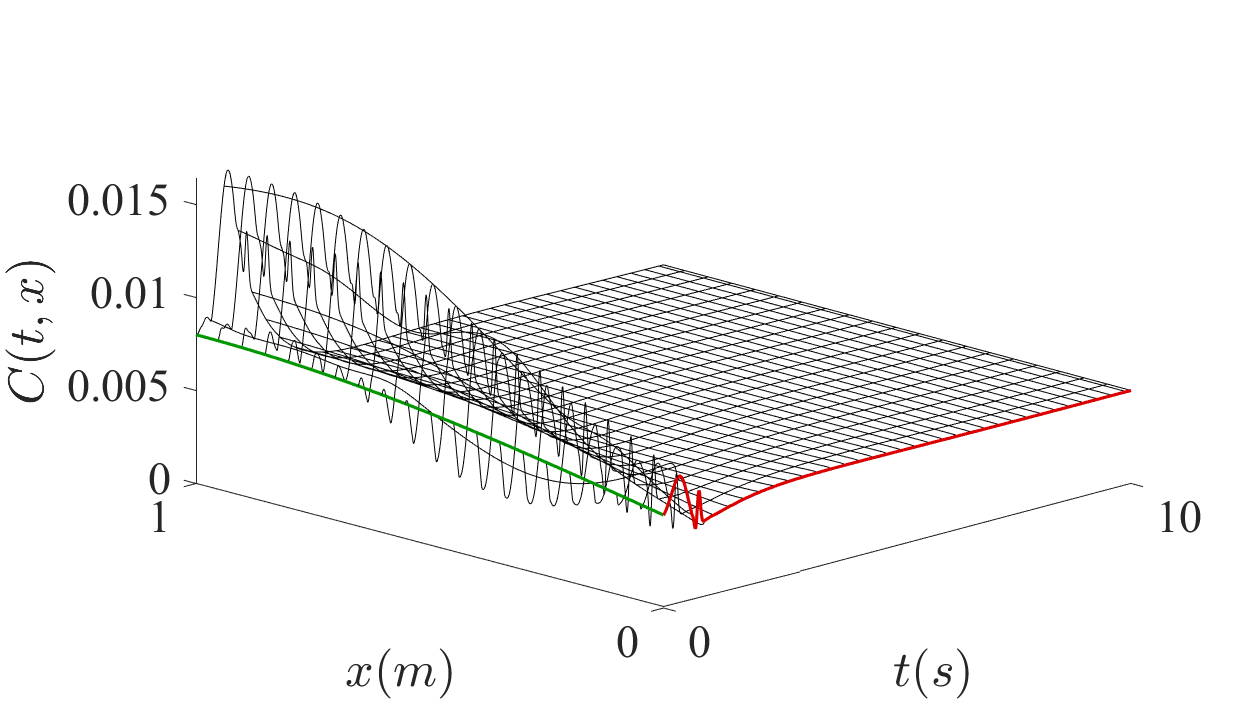}}
\caption{System variables under observer based feedback control}\label{fig:states}
\end{figure}

 The norm of the characteristics, $\|u_i\|_{\mathscr{L}^2}$ for $i=1,2,3$ and $\|w\|_{\mathscr{L}^2}$ are shown to converge to $0$ in Fig~\ref{fig:norm_cl}. As expected, the gate aperture of the underflow gate $U_L(t)$ is less than the water level at the gate, for the selected initial conditions and the system parameters. However, the open-loop system with  a setpoint control  defined as $U_L(t)=U_{eq}$ is  unstable in  $\mathscr{L}^2$ sense as depicted in Fig.~\ref{fig:norm_ol} for the same initial conditions given in \eqref{ini1}-\eqref{ini4}. Therefore, for this particular process, the  design of a feedback law  is crucial to maintain a healthy and sustainable water system.  The variation of the water level and the gate aperture  and that of the water velocity at the output gate, namely, the speed at which the water enters the downstream of the canal are shown in Fig.~\ref{fig:Control} and Fig.~\ref{fig:v_L}, respectively. The equilibrium values are reached around $t=8s$ under closed-loop control. The  evolution of the distributed states, $H(t,x), V(t,x), Z(t,x)$ and $C(t,x)$ with time $t$ and along the channel length $x$ are shown in Fig.~\ref{fig:states} to confirm the stability of the closed-loop system.

\begin{figure}
\centering
\subfloat[Entrainment rate\label{fig:E}]{
\includegraphics[width=0.45\columnwidth,trim={0 0 0 0.5cm},clip]{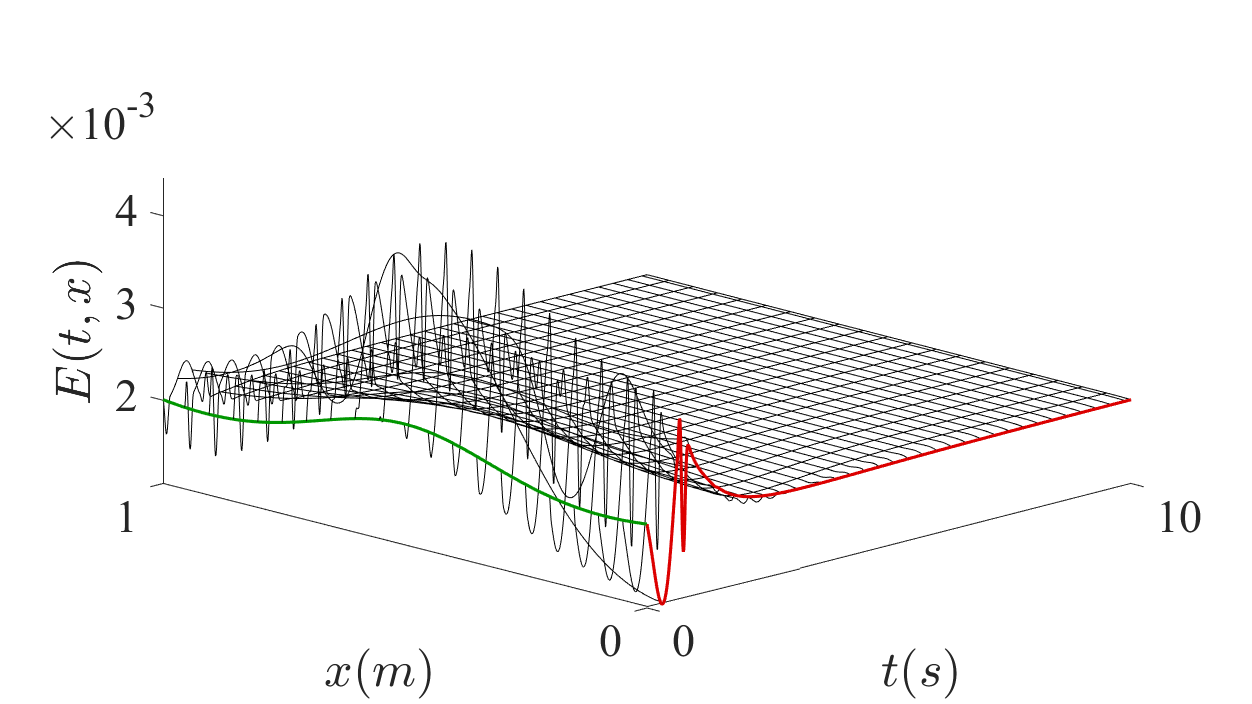}}
~
\subfloat[Deposition rate\label{fig:D}]{
\includegraphics[width=0.45\columnwidth,trim={0 0 0 0.5cm},clip]{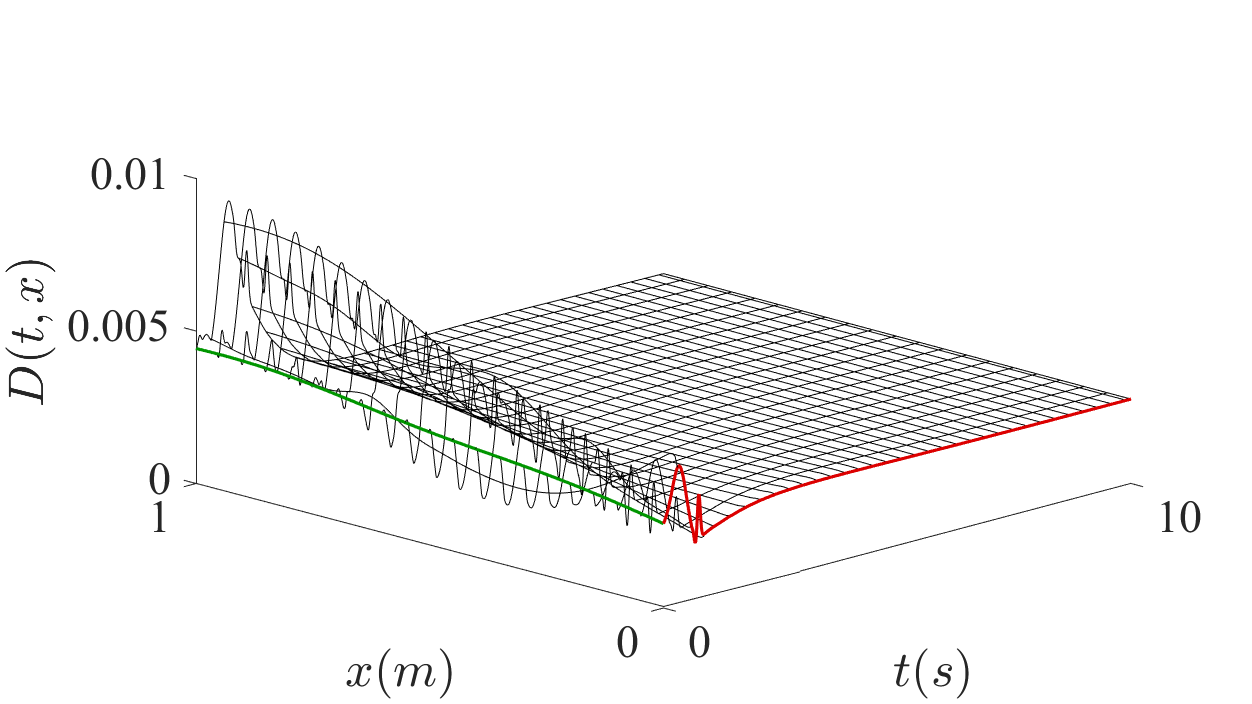}}
\caption{Sediment entrainment and deposition rates for the system}\label{fig:e and d}
\end{figure}

The dynamics of the spatially  distributed sediment entrainment and the deposition processes defined by the constitutive relations $E(V)$ and $D(V,C)$, respectively, are shown in Fig.~\ref{fig:E} and Fig.~\ref{fig:D} (we refer the reader to equations \eqref{erosion} and \eqref{deposition}). Initially, the sediment entrainment and deposition rates are notably high due to the selection of a low bed load equilibrium profile to stabilize. This choice implies an initial fast removal  of sediment particles from the bathymetry, as illustrated in Fig.~\ref{fig:z}. The rapid removal of the  excess sediment is illustrated in Fig.~\ref{fig:c}, where an immediate and substantial rise in the distributed suspended sediment concentration occurs right after the control signal kicks in. The significant increase in deposited sediment is notable because the sediment deposition process relies on both the suspended sediment concentration and the flow velocity, as denoted by equation \eqref{deposition}. Unlike the sediment deposition process, as indicated by equation \eqref{erosion}, the concentration of suspended sediment does not affect erosion. Instead, the rapid increase in eroded bedload is solely influenced by the flow velocity (a state variable) and other factors, including particle size, settling velocity, and the particle Reynolds number (refer to Appendix \ref{sec:app_1}).  In conclusion, during the transient phase, an initial high flow rate promotes high erosion, which  results in an increase in the suspended sediment concentration.  This, in turn, triggers an upsurge in the rate of sediment deposition all supported by a rapid increase in the flow velocity.

\begin{figure}
\centering
\subfloat[Mass flow rate of suspended sediment at the upstream and downstream gates \label{fig:sediment}]{
\includegraphics[width=0.45\columnwidth,clip]{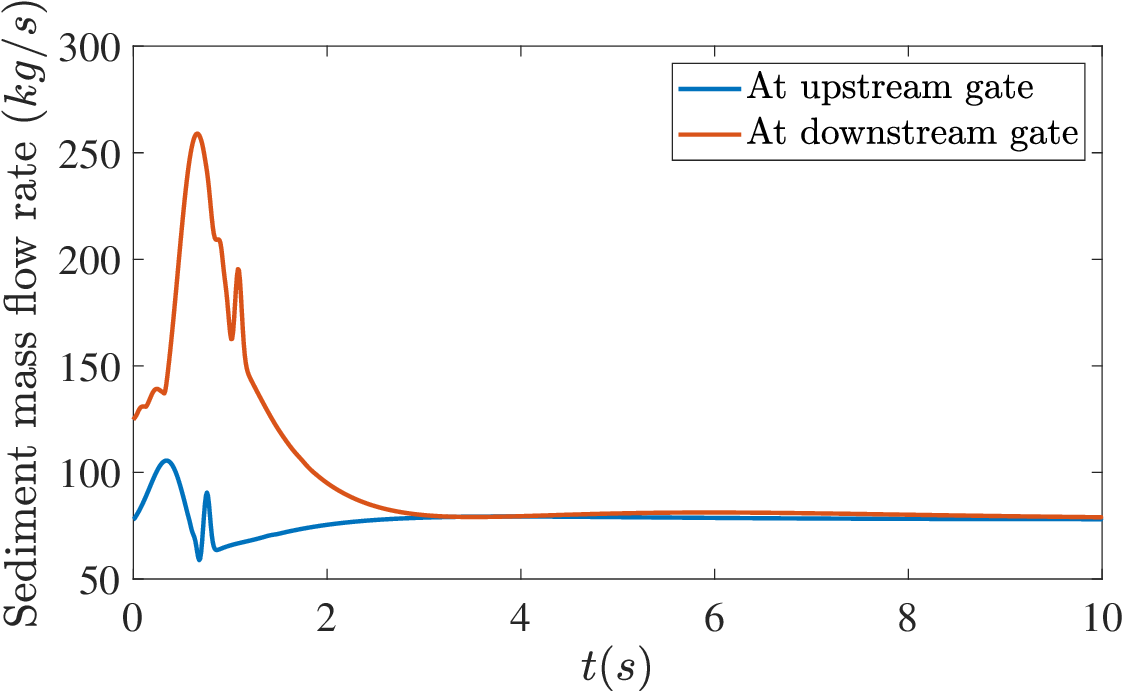}}
~
\subfloat[Variation of flushing effectiveness of sediment with time \label{fig:flushing_time}]{
\includegraphics[width=0.45\columnwidth,clip]{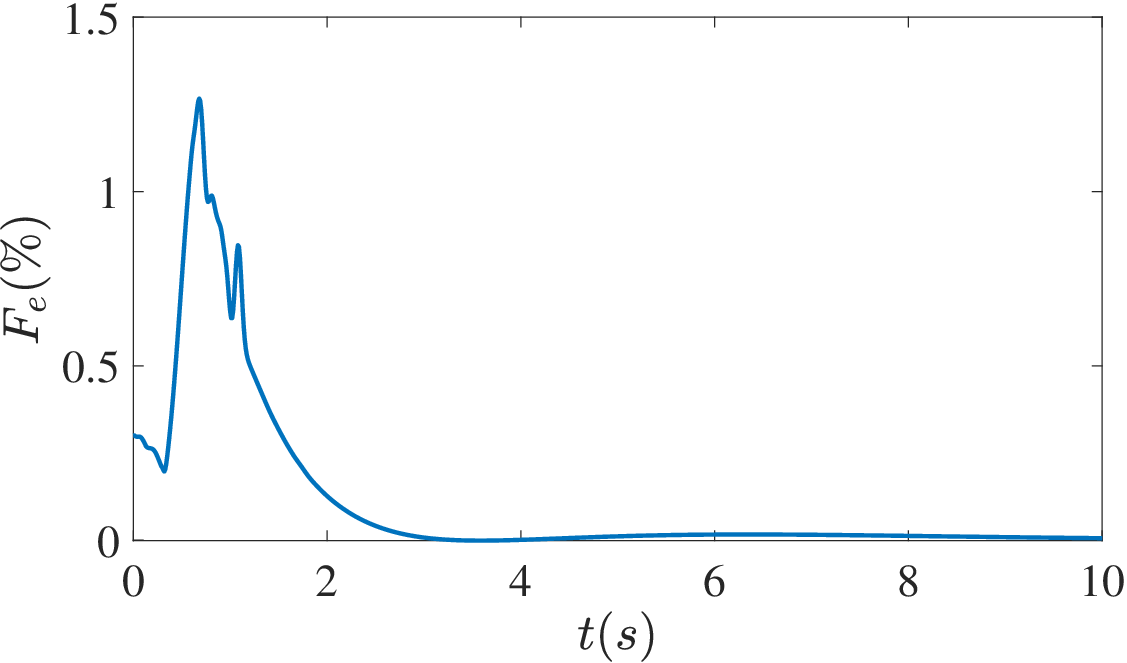}}
\caption{Suspended sediment motion and flushing effectiveness}\label{fig:suspended}
\end{figure}

Fig. \ref{fig:suspended} depicts the characteristics of suspended sediment dynamics. The mass flow rate of suspended sediment at the inlet and the outlet gates are shown in Fig.~\ref{fig:sediment}. Interestingly, the sediment mass flow rates exhibit an initial increase followed by a gradual decrease at the downstream gate. A similar trend  is observed  in Fig.~\ref{fig:flushing_time}, where the flushing effectiveness of the water channel, defined as $F_e=\frac{V_s}{V_w}$ is shown. This ratio represents the the net sediment volume removed ($V_s$) compared to the amount of released water ($V_w$) over a small time interval $\Delta t$ (Fig.~\ref{fig:flushing_time}). Higher sediment flow rates correspond to increased flushing effectiveness (\cite{jihn1996}). The flushing action of the controller is evident by the higher values of sediment mass flow rate at the downstream gate which converge to the upstream sediment flow rate as the excess sediment is flushed out. Overall, the simulation results indicate removal of sediment occurs at a relatively low suspended sediment ratio, reaching a maximum flushing effectiveness of $1.27\%$. The inefficient flushing of sediment is not considered a disadvantage but rather suggests the gradual removal of sediment as beneficial for preserving the river ecosystems. In addition, the modeling assumptions permit only small deviations from the equilibrium value of sediment concentration, thereby allowing for only low flushing effectiveness.
The exponential stability of the closed-loop system is confirmed by all simulation results with all states reaching the equilibrium set in  $8 s$.

\begin{figure}
    \centering
    \subfloat[Variation of $\mathscr{L}^2$-norm of $E(t,x)$ and $D(t,x)$ for the open-loop system\label{fig:EDu0}]{
    \includegraphics[width=0.45\columnwidth]{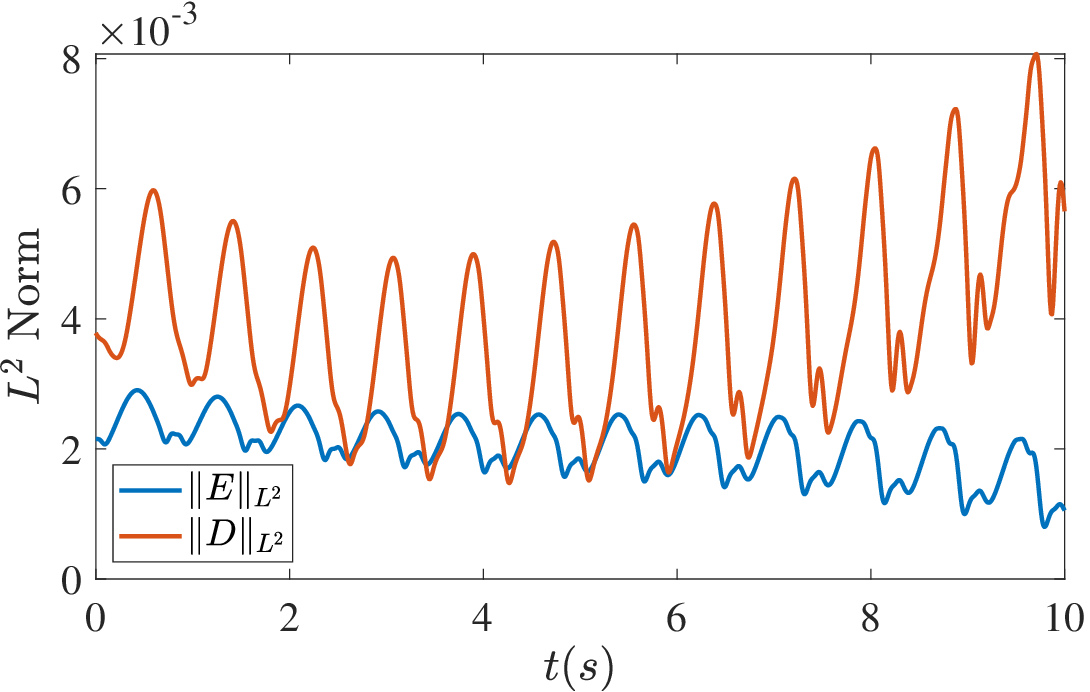}}
    ~
    \subfloat[Variation of flushing effectiveness of sediment with time for the open-loop system\label{fig:feu0}]{
    \includegraphics[width=0.45\columnwidth]{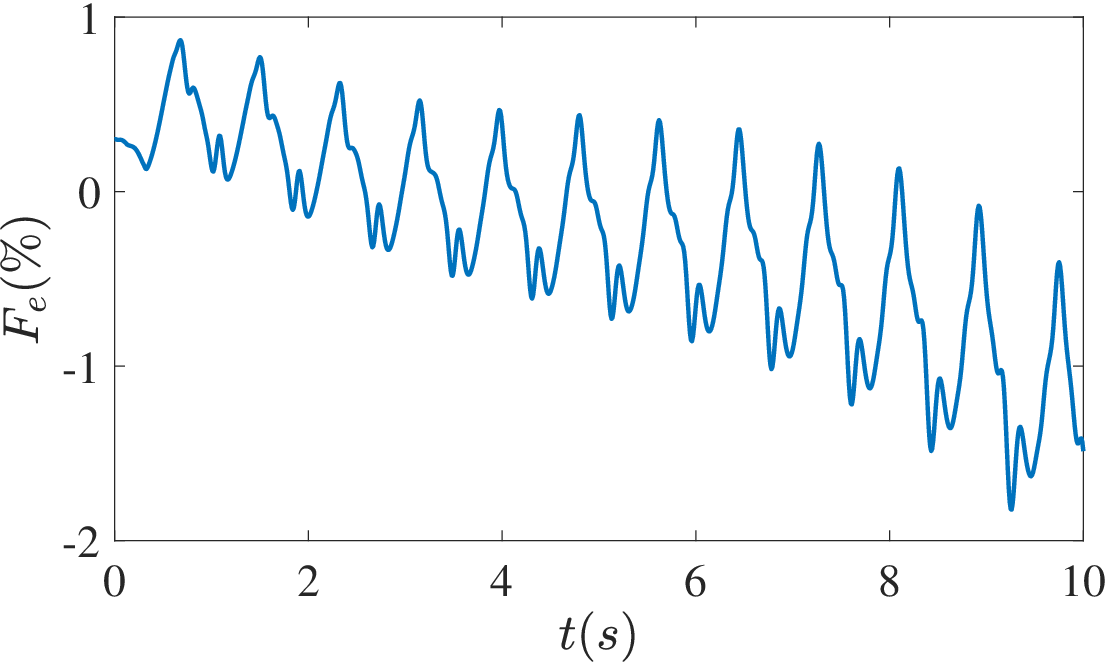}}
    \caption{Representation of sediment dynamics in the open-loop system}
\end{figure}

The norm of the characteristics of the open-loop system   in Fig.~\ref{fig:norm_ol}, shows  a drastic deviation from the equilibrium point which is incompatible with the assumptions of the linearized system. For the open-loop system, the time evolution of the $\mathscr{L}^2$-norm of the sediment entrainment and deposition rates as well as the flushing effectiveness,  is illustrated  in Fig.~\ref{fig:EDu0} and Fig.~\ref{fig:feu0}, respectively. The accumulation of sediment in the water canal under the open-loop system is apparent, driven by the higher sediment deposition rate in comparison to the sediment entrainment rate and the negative flushing effectiveness. Therefore, full automation of the downstream gate prevents degradation of the water canal since we obtain exponential stability of the system where the bathymetry gradually reaches the desired equilibrium point unlike in the open-loop case.

\section{Conclusion}\label{sec:conclusion}
 The main contribution  of this paper is the exponential stabilization of a new control-oriented model of coupled water-sediment flow in a channel  that takes the suspended sediment dynamics into consideration. We have shown that by measuring the state variables at the upstream gate, an anti-collocated boundary controller can be designed to exponentially stabilize the linearized model in  $\mathscr{L}^2$ sense.   The water height, the water velocity, the bathymetry and the suspended sediment concentration are stabilized to a desired setpoint. In addition, we have shown that without the automation of the downstream sluice gate, sediment may  accumulate over time, resulting in the degradation of the canal and its ecosystem. Considering that canals and river breaches are often exploited to store water in reservoirs that are in turn connected to water distribution networks, this accumulation of sediment over years could be viewed as a transgenerational issue. The proposed controller prevents the build up of sediment particles by regulating the dynamics of the system while maintaining the bathymerty at a desired level.

Overall, the modeling assumptions, namely, small deviations from the equilibrium and
the gradual variation of the bottom bed provide a starting point for analyzing sediment transport in open channels and serve as a foundation for further studies and applications. These assumptions may not always hold true in all situations, and more complex models might be required for specific cases where deviations from these assumptions are significant. Further developments will consider advance models that enables venting of turbidity via control design.

\begin{ack}
This work was funded by the NSF CAREER Award CMMI-2302030. The authors would like to thank Diwei Zhang for his contribution to the development of the control-oriented model with suspended sediment. 
\end{ack}

\appendix
\section{Empirical formulas} \label{sec:app_1}
The system parameters are defined using the morphodynamics model defined in \cite{garcia2008,Asselman2000}. Consider the parameters of the system given in table~\ref{tab:sys_parameters}.

\begin{table}[t]
\caption{Constant parameters of the system}
\begin{center}
\label{tab:sys_parameters}
\begin{tabular}{l >{\raggedright\arraybackslash}p{4cm} l}
& & \\ 
\hline
& Definition & Value\\
\hline
$g$&Acceleration due to gravity&$9.81m/s^2$\\
$C_f$&Bed friction coefficient&$0.002$\\
$k_G$&Discharge coefficient of sluice gate&0.6\\
$A$&Empirical constant&$1.3\times10^{-7}$\\
$c_1$&Empirical constant&$0.005$\\
$c_2$&Empirical constant&$1.3\times 10^{-11}$\\
$c_3$&Empirical constant&$2.75$\\
$r_1$&Empirical constant&$1$\\
$r_2$&Empirical constant&$31.5$\\
$r_3$&Empirical constant&$-1.46$\\
$\nu$&Kinematic viscosity at $20^\circ C$&$1\times 10^{-6} m^2/s$\\
$p^\prime$&Sediment porosity&$0.3$\\
$D_s$&Sediment size&$2.5\times 10^{-4} m$\\
$R$&Submerged specific gravity of sediment&$1.65$\\
$H_L$&Water height at downstream gate&$1 m$\\
\hline
\end{tabular}

\end{center}
\end{table}

The constant  $a$ in \eqref{PDE c},  the particle Reynolds number, $R_{ep}$ and  the sediment settling velocity $\nu_s$
are  defined as
\begin{align}
    a&=\frac{39C_f^{\frac{3}{2}}}{(1-p^\prime)gR},\\
    R_{ep}&=\frac{D_s\sqrt{RgD_s}}{\nu},\\
    \nu_s&=\frac{gRD_s^2}{18\nu},
\end{align}
respectively. The macro-parameters  $A_1$ and $A_2$ are given below
\begin{align}
    A_1=&A\left(\frac{\sqrt{C_f}}{\nu_s}R_{ep}^{0.6}\right)^5,\\
    A_2=&\frac{A_1}{0.3}.
\end{align}

The value of the constant bottom slope, $S_b$ is selected such that the relationship $gH_{eq}S_b=C_fV_{eq}^2$ is satisfied.

\section{Linearized model parameters} \label{sec:app_2}
The parameters defined in $\mathcal{B}$ are
\begin{align}
    \phi_{vh}=&\frac{C_fV_{eq}^2}{H_{eq}},\quad
    \phi_{vv}=\frac{-2C_fV_{eq}}{H_{eq}},\\
    \nonumber
    \phi_{zu}=&\frac{r_2r_3C_f^{\frac{r_3}{2}}}{1-p^\prime}\left(\frac{V_{eq}}{\nu_s}\right)^{r_3-1}\\
    &-\frac{5\nu_sA_1U_{eq}^4}{1-p^\prime}\left(\frac{1}{(1+A_2V_{eq}^5)^2}\right),\\
    \phi_{zc}=&\frac{v_s}{1-p^{\prime}}\left(r_1+r_2\left(\frac{\sqrt{C_f}}{\nu_s}\right)^{r_3}\right),\\
    \nonumber
    \phi_{ch}=&\frac{\nu_s}{H_{eq}^2}\bigg[C_{eq}\left(r_1+r_2\left(\frac{V_{eq}\sqrt{C_f}}{\nu_s}\right)^{r_3}\right)\\
    &-\frac{5A_1V_{eq}^4}{(1+A_2V_{eq}^5)^2}\bigg],\\
    \phi_{cu}=&\frac{5A_1V_{eq}^4}{H_{eq}(1+A_2V_{eq}^5)^2}-\frac{r_2r_3C_f^{\frac{r_3}{2}}}{H_{eq}}\left(\frac{V_{eq}}{\nu_s}\right)^{r_3-1},\\
    \phi_{cc}=&-\frac{\nu_s}{H_{eq}}\left(r_1+r_2\left(\frac{V_{eq}\sqrt{C_f}}{\nu_s}\right)^{r_3}\right).
\end{align}
The gate height at equilibrium, $U_{eq}$ can be calculated from \eqref{bc4} as
\begin{align}
    U_{eq}=\frac{Q_0}{k_G\sqrt{2g(H_{eq}+Z_{eq}-H_L)}}+Z_{eq}.\label{Ueq}
\end{align}
 The constants $\pi_{Lh},\pi_{Lv},\pi_{Lz},$ and $\pi_L$ are,
 \begin{align}
     \pi_{Lh}=&-\frac{V_{eq}(H_{eq}+2Z_{eq}-H_L)}{2(H_{eq}+Z_{eq}-H_L)},\\
     \pi_{Lv}=&-H_{eq},\\
     \nonumber
     \pi_{Lz}=&-k_G\sqrt{2g(H_{eq}+Z_{eq}-H_L)}\\&+\frac{Q_0}{2(H_{eq}+Z_{eq}-H_L)},\\
     \pi_L=&k_G\sqrt{2g(H_{eq}+Z_{eq}-H_L)}.
 \end{align}
\balance

\input{Conference.bbl}
\end{document}

%% file: figures/system.tex
\begin{tikzpicture}
\def\bccolor{RoyalBlue}
\def\concolor{BrickRed}
\def\coucolor{Sepia}

\def\len{0.3}
\node (u10) []{};
\node (u11) [right =\len\textwidth of u10, label={[xshift=-10 pt, yshift=0 pt]$u_1(x,t)$}] {};
\node (u20) [below =0.05\textwidth of u10] {};
\node (u21) [right =\len\textwidth of u20, label={[xshift=-10 pt, yshift=0 pt]$u_2(x,t)$}] {};
\node (u30) [below =0.05\textwidth of u20] {};
\node (u31) [right =\len\textwidth of u30, label={[xshift=-10 pt, yshift=0 pt]$u_3(x,t)$}] {};
\node (w0) [below =0.05\textwidth of u30] {};
\node (w1) [right =\len\textwidth of w0, label={[xshift=-10 pt, yshift=0 pt]$w(x,t)$}] {};
\node (y0) [left = 5pt of w0] {};
\node (y1) [left = 15pt of y0] {};
\node (U1) [right =5 pt of w1] {};
\node (U0) [right = 15pt of U1] {};
\node (cy) [below = 0.08\textwidth of y1] {};
\node (cU) [below = 0.08\textwidth of U0] {};
\node (c0) [right = 30pt of cy] {};
\node (obs) [rectangle, draw = \concolor, text = \concolor, minimum width = 2cm, minimum height = 1cm] at (c0) {Observer};
\node (c1) [left = 30pt of cU] {};
\node (con) [rectangle, draw = \concolor, text = \concolor, minimum width = 2cm, minimum height = 1cm] at (c1) {Controller};
\def\gap{0.03}
\node (u1u2) [right=\gap\textwidth of u10] {};
\node (u1u3) [right=\gap\textwidth of u1u2] {};
\node (u1w) [right=\gap\textwidth of u1u3] {};

\node (u2u1) [right=\gap\textwidth of u20] {};
\node (u2_1) [right=\gap\textwidth of u2u1] {}; 
\node (u2_2) [right=\gap\textwidth of u2_1] {}; 
\node (u2u3) [right=\gap\textwidth of u2_2] {};
\node (u2w) [right=\gap\textwidth of u2u3] {};

\node (u3_1) [right=\gap\textwidth of u30] {};
\node (u3u1) [right=\gap\textwidth of u3_1] {};
\node (u3_2) [right=\gap\textwidth of u3u1] {};
\node (u3u2) [right=\gap\textwidth of u3_2] {};
\node (u3_3) [right=\gap\textwidth of u3u2]{};
\node (u3w) [right=\gap\textwidth of u3_3] {};

\node (w_1) [right=\gap\textwidth of w0] {};
\node (w_2) [right=\gap\textwidth of w_1] {}; 
\node (wu1) [right=\gap\textwidth of w_2] {};
\node (w_3) [right=\gap\textwidth of wu1] {};
\node (wu2) [right=\gap\textwidth of w_3] {};
\node (wu3) [right=\gap\textwidth of wu2] {};

\node (x0) [below =5pt of w0, label={[xshift=0 pt, yshift=-20 pt]$x=0$}] {};
\node (x1) [below =5pt of w1, label={[xshift=0 pt, yshift=-20 pt]$x=L$}] {};
\node (x1a) [right =1pt of x1] {};
\draw [-stealth, line width=2pt]  (u10.center) to (u11.center);
\draw [-stealth, line width=2pt]  (u20.center) to (u21.center);
\draw [-stealth, line width=2pt] (u30.center) to (u31.center);
\draw [stealth-, line width=2pt] (w0.center) to (w1.center);
\draw [\bccolor,-stealth, line width=1pt,bend left=45] (w0.west) to node[pos=1.05,left]{$\delta_1$} (u10.west);
\draw [\bccolor,-stealth, line width=1pt,bend left=45] (w0.west) to node[pos=1.1,left]{$\delta_2$} (u20.west);
\draw [\bccolor,-stealth, line width=1pt,bend left=45] (w0.west) to node[pos=1.15,left]{$\delta_3$} (u30.west);
\draw [\bccolor,-stealth, line width=1pt,bend left=45]  (u11.east) to node[pos=0,right]{$\rho_1$} (w1.south east);
\draw [\bccolor,-stealth, line width=1pt,bend left=45]  (u21.east) to node[pos=0,right]{$\rho_2$} (w1.north east);
\draw [\concolor,stealth-, line width=1pt] (U1) to node[pos=0.8,above]{$U(t)$} (U1 -| cU);
\draw [\concolor,-, line width=1pt] (y0) to node[pos=0.8,above]{$w(t,0)$} (y0 -| cy);
\draw [\concolor,-, line width=1pt] (cy.center) to (cy |- y0);
\draw [\concolor,-,line width=1pt] (cU.center) to (cU |- U0);
\draw [\concolor,-stealth, line width=1pt] (cy.center) to (obs.west);
\draw [\concolor,-, line width=1pt] (cU.center) to (con.east);
\draw [\concolor,-stealth, line width=1pt] (obs.east) to node[pos=0.5, above] {$\hat{u}(t,x),\hat{w}(t,x)$} (con.west);
\draw [\coucolor,stealth-stealth,dashed, line width=1pt] (u1u2.south) to node[pos=0.2,left]{$\sigma_{12}$} node[pos=0.8,left]{$\sigma_{21}$} (u2u1.north);
\draw [\coucolor,stealth-stealth,dashed, line width=1pt] (u1u3.south) to node[pos=0.18,left]{$\sigma_{13}$} node[pos=0.9,left]{$\sigma_{31}$} (u3u1.north);
\draw [\coucolor,stealth-stealth,dashed, line width=1pt] (u1w.south) to node[pos=0.1,left]{$\alpha_1$} node[pos=0.9,left]{$\theta_1$} (wu1.north);
\draw [\coucolor,stealth-stealth,dashed, line width=1pt] (u2u3.south) to node[pos=0.2,left]{$\sigma_{23}$} node[pos=0.9,left]{$\sigma_{32}$} (u3u2.north);
\draw [\coucolor,stealth-stealth,dashed, line width=1pt] (u2w.south) to node[pos=0.2,left]{$\alpha_2$} node[pos=0.9,left]{$\theta_2$} (wu2.north);
\draw [\coucolor,stealth-stealth,dashed, line width=1pt] (u3w.south) to node[pos=0.2,left]{$\alpha_3$} node[pos=0.9,left]{$\theta_3$} (wu3.north);
\draw [-stealth, line width=0.8pt] (x0.west) to (x1a.east);
\draw [line width=0.5pt] (x0.north) -- (x0.south);
\draw [line width=0.5pt] (x1.north) -- (x1.south);
\end{tikzpicture}